\documentclass[12pt]{amsart}

\usepackage[latin1]{inputenc}
\usepackage[english]{babel}
\usepackage{amsthm}
\usepackage{amsmath}
\usepackage{amsxtra}
\usepackage{mathrsfs}
\usepackage{amssymb}
\usepackage[dvips]{graphicx}
\usepackage{graphicx}

\input xy
\xyoption{all}

\newtheorem{teo}{Theorem}[section]
\newtheorem{lemma}[teo]{Lemma}

\newtheorem{defi}[teo]{Definition}
\newtheorem{coro}[teo]{Corollary}
\newtheorem{fatto}[teo]{Fact}

\newtheorem{cla}[teo]{Claim}
\newtheorem{prop}[teo]{Proposition}

\theoremstyle{remark}
\newtheorem{remark}[teo]{Remark}

\begin{document}

\newcommand{\ran}{{\rm ran}}
\newcommand{\cof}{{\rm cof}}
\newcommand{\dom}{{\rm dom}}
\newcommand{\I}{I}
\newcommand{\N}{\mathcal{N}}
\newcommand{\up}{\upharpoonright}
\newcommand{\urltilde}{\kern -.15em\lower .7ex\hbox{~}\kern .04em}

\newcommand\M{\mathcal{M}}
\newcommand\Acal{\mathscr{A}}
\newcommand\Bcal{\mathscr{B}}
\newcommand\Ecal{\mathscr{E}}
\newcommand\Fcal{\mathscr{F}}
\newcommand\Hcal{\mathscr{H}}
\newcommand\Ical{\mathscr{I}}
\newcommand\Ncal{\mathscr{N}}
\newcommand\Mcal{\mathscr{M}}
\newcommand\Pcal{\mathscr{P}}
\newcommand\Qcal{\mathscr{Q}}
\newcommand\Rcal{\mathscr{R}}
\newcommand\Tcal{\mathscr{T}}
\newcommand\Ucal{\mathscr{U}}
\newcommand\Zcal{\mathscr{Z}}
\newcommand\Rbb{\mathbb{R}}
\newcommand\Nfrak{\mathfrak{N}}
\newcommand\Pfrak{\mathfrak{P}}
\newcommand\restrict{\restriction}

\newcommand{\diff}{\operatorname{\mathrm diff}}
\newcommand{\Ht}{\operatorname{\mathrm ht}}
\newcommand{\lev}{\operatorname{\mathrm lev}}
\newcommand{\meet}{\wedge}
\newcommand\triord{\triangleleft}
\newcommand{\Th}{{}^{\mathrm{th}}}
\newcommand{\St}{{}^{\mathrm{st}}}
\newcommand\axiom{\mathrm}
\newcommand\MA{\axiom{MA}}
\newcommand\MM{\axiom{MM}}
\newcommand\PFA{\axiom{PFA}}
\newcommand\BPFA{\axiom{BPFA}}
\newcommand\MRP{\axiom{MRP}}
\newcommand\SRP{\axiom{SRP}}
\newcommand\ZFC{\axiom{ZFC}}
\newcommand\CAT{\axiom{CAT}}
\newcommand\CH{\axiom{CH}}
\newcommand{\<}{\langle}
\renewcommand{\>}{\rangle}
\newcommand\mand{\textrm{ and }}
\renewcommand{\diamond}{\diamondsuit}
\newcommand{\forces}{\Vdash}

\newcommand{\cf}{{\mbox{cof}}}
\newcommand{\height}{\Ht}
\newcommand\NS{\mathrm{NS}}

\title[Proper forcing remastered]{Proper forcing remastered}
\author[Boban Veli\v{c}kovi\'{c}]{Boban Veli\v{c}kovi\'{c}}
\email{boban@logique.jussieu.fr}
\urladdr{http://www.logique.jussieu.fr/\urltilde boban}
\author{Giorgio Venturi}
\email{gio.venturi@gmail.com}
\address{Equipe de Logique,
Universit\'e Paris 7 Diderot,
site Chevaleret 75205, Paris Cedex 13, France}

\keywords{forcing, proper, side conditions, club}

\maketitle

\begin{abstract}{ We present the method introduced by Neeman of generalized side conditions
with two types of models. We then discuss some applications: a variation of
the Friedman-Mitchell poset for adding a club with finite conditions, the consistency
of the existence of an $\omega_2$ increasing chain in $(\omega_1^{\omega_1},<_{\rm fin})$, originally
proved by Koszmider, and the existence of a thin very tall superatomic Boolean algebra,
originally proved by Baumgartner-Shelah. We expect that the present method
will have many more applications.}
\end{abstract}

\section*{Introduction}

We present a generalization of the method of model as side conditions.
Generally speaking a poset that uses models as side conditions is a notion
of forcing whose elements are pairs, consisting of a working part which is
some partial information about the object we wish to add and   a finite
$\in$-chain of countable elementary substructures of $H(\theta)$, for some cardinal
$\theta$ i.e. the structure consisting of sets whose transitive closure has
cardinality less than $\theta$. The models in the side condition are used
to control the extension of the working part. This is crucial
in showing some general property of the forcing such as properness.

The generalization we now present amounts to allowing also certain uncountable models
in the side conditions. This is used to show that the forcing preserves both $\aleph_1$ and $\aleph_2$.
This approach was introduced by Neeman \cite{Neeman} who used it to give an alternative
proof of the consistency of PFA and also to obtain generalizations of PFA
to higher cardinals. In \S 1 we present the two-type poset of pure side conditions from \cite{Neeman},
in the case of countable models and approachable models of size $\omega_1$,
and work out the details of some of its main properties that were mentioned in \cite{Neeman}.
The remainder of the paper is devoted to applications. We will be primarily interested
in adding certain combinatorial objects of size $\aleph_2$.
These results were known by other methods but we believe that the present method
is more efficient and will have other applications.  In \S 2 we present a version of
the forcing for adding a club in $\omega_2$ with finite conditions, preserving $\omega_1$
and $\omega_2$. This fact has been shown to be consistent with ZFC independently
by Friedman (\cite{Fried}) and Mitchell (\cite{Mitch}) using  more complicated
notions of forcing. In \S 3 we show how to add a  chain of
length $\omega_2$ in the structure $(\omega_1^{\omega_1},<_{\rm fin})$.
This result is originally due to Koszmider \cite{Koszmider}.
Finally, in \S 4 we give another proof of a result of Baumgartner and Shelah
\cite{BaumShelah} by using side condition forcing to add a thin very tall superatomic
Boolean algebra.

This paper represents the content of the lectures of the first author given
at the Appalachian Set Theory Workshop in Chicago, October 15 2011.
The second author has helped in the preparation of these notes.

\section{The forcing $\mathbb{M}$}

In this section we present the forcing consisting of pure side conditions.
Our presentation follows \cite{Neeman}, but we only consider side conditions
consisting of models which are either countable or of size $\aleph_1$.
We consider the structure $(H(\aleph_2),\in, \unlhd )$ equipped with a fixed
well-ordering $\unlhd$. In this way we have definable Skolem functions,
so if $M$ and $N$ are elementary submodels of $H(\aleph_2)$ then so is $M\cap N$.

\begin{defi}
Let $P$ an elementary submodel of $H(\aleph_2)$ of size $\aleph_1$. We say that $P$ is
\emph{internally approachable} if it can be written as the union of an increasing
continuous $\in$-chain  $\<P_{\xi}: \xi < \omega_1\>$ of countable
elementary submodels of $H(\aleph_2)$ such that $\<P_{\xi}: \xi < \eta\> \in P_{\eta +1}$,
for every ordinal $\eta < \omega_1$.
\end{defi}

If $P$ is internally approachable of size $\aleph_1$ we let $\vec{P}$
denote the least $\unlhd$-chain witnessing this fact and we write $P_\xi$ for the $\xi$-th
element of this chain. Note also that in this case  $\omega_1\subseteq P$.

\begin{defi}
We let $\mathcal{E}_0^{2}$ denote the collection of all countable elementary submodels of $H(\aleph_2)$
and $\mathcal{E}_1^{2}$ the collection of all internally approachable elementary submodels of
$H(\aleph_2)$ of size $\aleph_1$. We let $\mathcal E^2=\mathcal E ^2_0\cup \mathcal E_1^2$.
\end{defi}

The following fact is well known.

\begin{fatto}
The set $\mathcal{E}_1^{2}$ is stationary in $[H(\aleph_2)]^{\aleph_1}$.
\qed \end{fatto}

We are now ready to define the forcing notion $\mathbb M$ consisting of pure side conditions.

\begin{defi}
The forcing notion $\mathbb{M}$ consists of finite $\in$-chains $p = \mathcal{M}_p$ of models in
$\mathcal{E}^2$  closed under intersection.
The order on $\mathbb{M}$ is reverse inclusion, i.e. $q \leq p$ if $\mathcal{M}_p \subseteq \mathcal{M}_q$.
\end{defi}

\medskip

Suppose $M$ and $N$ are elements of $\mathcal{E}^{2}$ with $M\in N$.
If $|M|\leq |N|$ then $M\subseteq N$. However, if $M$ is of size $\aleph_1$ and $N$
is countable then the $\unlhd$-least chain $\vec{M}$ witnessing that $M$ is internally approachable
belongs to $N$ and so $M\cap N= M_{\delta_N}$, where $\delta_N =N\cap \omega_1$ and $M_{\delta_N}$
is the $\delta_N$-th member of $\vec{M}$.

\medskip

We can split every condition in $\mathbb{M}$ in two parts:
the models of size $\aleph_0$ and the models of size $\aleph_1$.

\begin{defi} For $p \in \mathbb M$ let
 $\pi_0(p)=p\cap \mathcal{E}_0^{2}$ and $\pi_1(p)=p\cap \mathcal{E}_1^{2}$.
\end{defi}

Let us see some structural property of the elements of $\mathbb{M}$.
First, let $\in ^*$ be the transitive closure of the $\in$ relation, i.e.
$x\in ^* y$ if $x\in {\rm tcl}(y)$. Clearly, if $p \in \mathbb M$ then
$\in ^*$ is a total ordering on $\mathcal M_p$. Given
$M,N \in \mathcal M_p \cup \{ \emptyset, H(\aleph_2) \}$ with $M\in ^* N$ let
$$
(M,N)_p = \{ P\in \mathcal M_p: M\in ^* P \in ^* N\}.
$$
We let $(M,N]_p=(M,N)_p\cup \{ N\}$, $[M,N)_p =(M,N)_p \cup \{ M\}$ and
$[M,N]_p=(M,N)_p \cup \{ M,N\}$.
Given a condition $p\in \mathbb M$ and $M\in p$ we let
$p\restriction M$ denote the restriction of $p$ to $M$, i.e.
$\mathcal M_p\cap M$.

\begin{fatto}\label{unctble-restriction}
Suppose $p \in \mathbb{M}$ and $N\in \pi_1(p)$. Then $\mathcal M_p \cap N = (\emptyset, N)_p$.
Therefore, $p \cap N \in \mathbb M$.
\qed
\end{fatto}

\begin{fatto}\label{ctble-restriction}
Suppose $p \in \mathbb{M}$ and $M \in \pi_0(p)$. Then
$$
\mathcal M_p \cap M = \mathcal M_p \setminus \bigcup \{ [M\cap N,N)_p : N\in (\pi_1(p)\cap M)\cup \{ H(\aleph_2)\} \}.
$$
Therefore, $p\cap  M \in \mathbb M$.
\qed
\end{fatto}

The next lemma will be used in the proof of properness of $\mathbb{M}$.

\begin{lemma}\label{countable}
 Suppose $M \in \mathcal E^2$ and $p\in  \mathbb M\cap M$.
 Then  there is a new condition  $p^M$, which is the smallest element of $\mathbb{M}$
 extending $p$ and containing $M$ as an element.
\end{lemma}
\begin{proof}
If $M\in \mathcal E_1^2$ we can simply let
$$
p^M=\mathcal M_{p^M}=\mathcal M_p \cup \{M\}.
$$
If $M\in \mathcal E_0^2$ we close $\M_{p} \cup \{M\}$ under intersections and
show that it is still an $\in$-chain. First of all notice that,
since $p$ is finite and belongs to $M$, we have $\M_p \subseteq M$.
For this reason if $P \in \pi_0(p)$, then $P \cap M = P$. On the other hand,
if $P \in \pi_1(p)$, by the internal approachability of $P$ and the fact
that $P\in M$ we have that $P \cap M \in  P$.
Now, if $N \in P$ is the $\in ^*$-greatest element of $\M_p$ below $P$,
then $N \in P \cap M$, since $\M_p \subseteq M$.
Finally the $\in ^*$-greatest element of $\M_p$ belongs to $M$,
since $\M_p$ does.
\end{proof}

\medskip

Let $\mathcal P$ be a forcing notion. We say that a set $M$ is {\em adequate} for
$\mathcal P$ if for every $p,q\in M\cap \mathcal P$ if $p$ and $q$ are compatible
then there is $r\in \mathcal P \cap M$ such that $r\leq p,q$.
Note that we do not require that $\mathcal P$ belongs to $M$.
In the forcing notions we consider if two conditions $p$ and $q$ are compatible
then this will be witnessed by a condition $r$ which is $\Sigma_0$-definable from $p$ and $q$.
Thus, all elements of $\mathcal E^2$ will be adequate for the appropriate forcing notions.

\begin{defi}
Suppose $\mathcal P$ is a forcing notion and $M$ is adequate for $\mathcal P$.
We say that a condition $p$ is {\em $(M,\mathcal P)$-strongly generic} if $p$ forces that $\dot{G} \cap M$
is a $V$-generic subset of $\mathcal P \cap M$, where $\dot{G}$ is the canonical
name for the $V$-generic filter over $\mathcal P$.
\end{defi}

In order to check that a condition is strongly generic over a set $M$
we can use the following characterization, see \cite{Mitch2} for a proof.

\begin{fatto}\label{char}
Suppose $\mathcal P$ a notion of forcing and $M$ is adequate for $\mathcal P$.
A condition $p$ is $(M,\mathcal P)$-strongly generic if and only if for every $r \leq p$ in $\mathcal P$
there is a condition $r | M \in \mathcal P \cap M$ such that any condition
$q \leq r | M$ in $M$ is compatible with $r$.
\end{fatto}

\begin{defi}
Suppose $\mathcal P$ is a forcing notion
and $\mathcal S$ is a collection of sets adequate for $\mathcal P$.
We  say that $\mathcal P$ is {\em $\mathcal S$-strongly proper}, if for every $M \in \mathcal S$,
every condition $p \in \mathcal P\cap  M$ can be extended to an $(M,\mathcal P)$-strongly
generic condition $q$.
\end{defi}

Our goal is to show that $\mathbb M$ is $\mathcal E^2$-strongly proper. We will need
the following.

\begin{lemma}\label{amalgam} Suppose $r\in \mathbb M$ and $M\in \mathcal M_r$. Let $q \in M$
be such that $q\leq r\cap M$. Then $q$ and $r$ are compatible.
\end{lemma}



\begin{proof}
If $M$ is uncountable then one can easily check that
$\mathcal M_s = \mathcal M_q \cup \mathcal M_r$ is an $\in$-chain which
is closed under intersection. Therefore $s=\mathcal M_s$ is a common extension
of $q$ and $r$.
Suppose now $M$ is countable.
We first check that $\M_{q} \cup \M_r$ is an $\in$-chain,
then we close this chain under intersections and show that the resulting set
is still an $\in$-chain.

\begin{cla}\label{prop}
The set $\M_{q} \cup \M_r$ is an $\in$-chain.
\end{cla}
\begin{proof}
Note that any model of $\mathcal M_r \setminus M$ is either in
$[M,H(\aleph_2))_r$ or belongs to an interval of the form
$[N\cap M,N)_r$, for some $N\in \pi_1(r\restriction M)$. Consider one such interval
$[N\cap M,N)_r$. Since $N\in r\restriction M$ and $q\leq r\restriction M$
we have that $N\in \mathcal M_q$. The models in $\mathcal M_r \cap [N\cap M,N)_r$
are an $\in$-chain. The least model on this chain is $N\cap M$
and the last one belongs to $N$.
Consider the $\in ^*$-largest model $P$ of $\mathcal M_q$ below $N$.
Since $q\in M$ we have that $P\in M$. Moreover, since $\mathcal M_q$
is an $\in$-chain we have that also $P\in N$, therefore $P\in N\cap M$.
Similarly, the least model of $\mathcal M_r$ in $[M,H(\aleph_2))_r$
is $M$ and it contains the top model of $\mathcal M_q$.
Therefore, $\mathcal M_q \cup \mathcal M_r$ is an $\in$-chain.
\end{proof}

We now close $\M_q \cup \M_r$ under intersections and check that it is still an $\in$-chain.
We let $Q \in \mathcal M_q \setminus \mathcal M_r$ and consider models of the form $Q\cap R$, for
$R\in \mathcal M_r$.
\medskip

\noindent \emph{Case 1}: $Q \in \pi_0(q)$. We show by $\in ^*$-induction
on $R$ that $Q\cap R$ is already on the chain $\mathcal M_q$.
Since $Q \in M$ and $Q$ is countable we have that $Q\subseteq M$.
Therefore, $Q\cap R = Q\cap (R\cap M)$.
We know that $R,M\in \mathcal M_r$ and $\mathcal M_r$ is closed under intersections,
so $R\cap M \in \mathcal M_r$. By replacing $R$ by $R\cap M$ we may assume that
$R$ is countable and below $M$ in $\mathcal M_r$. If $R\in M$ then $R\in \mathcal M_q$
and $\mathcal M_q$ is closed under intersection, so $Q\cap R \in \mathcal M_q$.
If $R \in \mathcal M_r \setminus M$ then it belongs to an interval of
the form $[N\cap M,N)_r$, for some $N\in \pi_1(r\restriction M)$.
Since $N$ is uncountable and $R \in ^*N$ it follows that $R\subseteq N$.
If there is no uncountable model in the interval $[N\cap M,R)_r$ then
we have that $N\cap M\subseteq R \subseteq N$.
It follows that
$$
Q\cap (N\cap M) \subseteq Q\cap R \subseteq Q\cap N.
$$
However, $Q$ is a subset of $M$ and so $Q\cap (N\cap M)= Q\cap N$.
Therefore, $Q\cap R=Q\cap N$ and since $Q, N \in \mathcal M_q$
we have again that $Q\cap N \in \mathcal M_q$.
Now, suppose there is an uncountable model in $[N\cap M,R)_r$ and let $S$
be the $\in^*$-largest such model.
Since all the models in the interval $(S,R)_r$ are countable we have that
$S\in R$. On the other hand, $S$ is uncountable and above $N\cap M$ in $\mathcal M_r$.
It follows that $N\cap M \subseteq S$.
Now, consider the model $R^*=R\cap S$. It is below $S$ in $\mathcal M_r$.
We claim that $Q\cap R=Q\cap R^*$.
To see this note that, since $Q\subseteq M$ and $R\subseteq N$, we have
$$
Q\cap R\subseteq Q\cap (N \cap M) \subseteq Q\cap S.
$$
Therefore, $Q\cap R^*= Q\cap (R\cap S)=Q\cap R$.
Since $R^*$ is below $R$ in $\mathcal M_r$, by the inductive assumption,
we have that $Q\cap R^* \in \mathcal M_q$.

\medskip

\noindent \emph{Case 2}: $Q\in \pi_1(q)$. We first show that the largest element of
$\mathcal M_q \cup \mathcal M_r$ below $Q$ is in $\mathcal M_q$.
To see this note that by Fact \ref{ctble-restriction} any model, say $S$,
in $\mathcal M_r \setminus M$ which is below $M$ under $\in ^*$
belongs to an interval of the form $[N\cap M,N)_r$, for some
$N\in \pi_1(r\restriction M)$. By our assumption, $Q \in \mathcal M_q \setminus \mathcal M_r$
so $N$ is distinct from $Q$. Since $N,Q \in \mathcal M_q$ and they are
both uncountable it follows that either $Q\in N$ or $N\in Q$.
In the first case, $Q\in N\cap M$, i.e. $Q$ is $\in ^*$-below $S$.
In the second case, $S\in ^* N\in ^* Q$ and $N\in M$.

We now consider models of the form $Q\cap R$, for $R\in \mathcal M_r$.
If $R$ is uncountable then either $Q\subseteq R$ or $R\subseteq Q$
so $Q\cap R$ is in $\mathcal M_q \cup \mathcal M_r$.
If $R$ is countable and below  $Q$ on the chain $\mathcal M_q \cup \mathcal M_r$
then $R\subseteq Q$, so $Q\cap R =R$.
If $R\in \mathcal M_r \cap M$ then $R\in \mathcal M_q$ and since
$\mathcal M_q$ is closed under intersections we have that $Q\cap R\in \mathcal M_q$.
So, suppose $R\in \pi_0(r)\setminus M$. By Fact \ref{ctble-restriction}
we know that $R$ is either in $[M,H(\aleph_2))_r$ or in $[N\cap M,N)_r$,
for some $N\in \pi_1(r\restriction M)$. We show by $\in ^*$-induction
that $Q\cap R$ is either in $\mathcal M_q\cup \mathcal M_r$ or
is equal to $Q_{\delta _R}$ and moreover  $\delta_R \geq \delta_M$.
Consider the case $R\in [M,H(\aleph_2))_r$. If there is no uncountable
$S$ in the interval $(M,R)_r$ then $M\subseteq R$. Therefore,
$Q\in R$ and $\delta_R \geq \delta_M$. Since $Q\in R$ then
$Q\cap R=Q_{\delta_R}$. If there is an uncountable model in
the interval $(M,R)_r$ let $S$ be the largest such model.
Since $Q$ is below $S$ in the $\mathcal M_q \cup \mathcal M_r$ chain
we have $Q\subseteq S$, so if we let $R^*=R\cap S$, then
$Q\cap R^* = Q\cap R$, and moreover $\delta_{R^*}=\delta_R$.
Therefore, we can use the inductive hypothesis for $R^*$.
The case when $R$ belongs to an interval of the form
$[N\cap M,N)_r$, for some $N\in \pi_1(r\restriction M)\cup \{ H(\aleph_2)\}$ is treated
in the same way.

The upshot of all of this is that when we close $\mathcal M_q\cup \mathcal M_r$
under intersections the only new models we add are of the form
$Q_{\xi}$, for $Q\in \pi_1(\mathcal M_q \setminus \mathcal M_r)$, and finitely many
countable ordinals $\xi\geq \delta_M$. These models form an $\in$-chain,
say $C_Q$.
In particular, the case $R = M$ falls under the last case of the previous paragraph,
therefore $Q_{\delta_M} = Q \cap M$ is the $\in^*$-least member of $C_Q$.
Moreover, if $Q'$ is the predecessor of $Q$ in $\mathcal M_q \cup \mathcal M_r$,
then $Q'$ belongs to both $Q$ and $M$ and hence it belongs to $Q_{\delta_M}$.
The largest member of $C_Q$ is a member of $Q$ since it is of the form
$Q_\xi$, for some countable $\xi$.
Thus, adding all these chains to $\mathcal M_q \cup \mathcal M_r$ we
preserve the fact that we have an $\in$-chain.
\end{proof}

As an immediate consequence of Lemma \ref{amalgam}
we have the following.

\begin{teo}\label{strongly-proper} $\mathbb M$ is $\mathcal E^2$-strongly proper.
\end{teo}

\begin{proof}
Suppose $M\in \mathcal E^2$ and $p\in M\cap \mathbb M$.
We shall show that $p^M$ is $(M,\mathbb M)$-strongly generic.
To see this we for every condition $r\leq p^M$  we  have to define
a condition $r|M\in \mathbb M \cap M$ such that
for every $q\in \mathbb M \cap M$
if $q\leq r|M$ then $q$ and $r$ are compatible.
If we let $r|M$ simply be $r\cap M$ this is precisely
the statement of  Lemma \ref{amalgam}.
\end{proof}

\begin{coro}
The forcing $\mathbb{M}$ is proper and preserves $\omega_2$.
\qed
\end{coro}

\section{Adding a club in $\omega_2$ with finite conditions}

We now present a version of the Friedman-Mitchell (see \cite{Fried} and \cite{Mitch})
forcing for adding a club to $\omega_2$ with finite conditions. This will be achieved
by adding a working part to the side conditions.

\begin{defi}\label{M_2definition}
Let $\mathbb{M}_2$ be the forcing notion whose elements are triples $p = (F_p, A_p, \mathcal{M}_p)$, where
$F_p \in [\omega_2]^{<\omega}$, $A_p$ is a finite collection of intervals of the
form $(\alpha, \beta]$, for some $\alpha, \beta <\omega_2$, $\mathcal{M}_p \in \mathbb{M}$, and
\begin{enumerate}
\item $F_p \cap \bigcup A_p=\emptyset$,
\item if $M \in \M_p$ and $I \in A_p$, then either $I \in M$ or $I \cap M = \varnothing$.
\end{enumerate}

\noindent The order on $\mathbb{M}_2$ is coordinatewise reverse inclusion, i.e. $q \leq p$ if $F_p\subseteq F_q$,
$A_p\subseteq A_q$ and $\mathcal M_p \subseteq \mathcal M_q$.
\end{defi}

The information carried by a condition $p$ is the following.
The points of $F_p$ are going to be in the generic club, and the intervals in $A_p$ are
a partial description of the complement of that club.
The side conditions are there to ensure that the forcing is $\mathcal E^2$-strongly proper.
It should be pointed our that a condition $r$ may force some
ordinals to be in the generic club even though they are not
explicitly in $F_r$. The reason is that we may not be able
to exclude them by intervals which satisfy conditions (1) and (2)
of Definition~\ref{M_2definition}.

\begin{fatto}\label{exclude-sup}
If $p \in \mathbb M_2$ and $M\in \mathcal M_p$ then $\sup(M\cap \omega_2)\notin \bigcup A_p$.
\end{fatto}

\begin{proof} Any interval $I$ which contains $\sup(M\cap \omega_2)$
would have to intersect $M$ without being an element of $M$. This contradicts
condition (2) of Definition~\ref{M_2definition}.
\end{proof}

\begin{fatto}\label{exclude-int}
Suppose $p\in \mathbb M_2$, $M\in \mathcal M_p$ and $\gamma \in F_p$. Then
$$
\min(M\setminus \gamma), \sup(M\cap \gamma)\notin \bigcup A_p.
$$
\end{fatto}

\begin{proof} Suppose $\gamma \in F_p$ and let $I\in A_p$.
Then $I$ is of the form $(\alpha,\beta]$, for some ordinals $\alpha,\beta <\omega_2$.
Since $p$ is a condition we know that $\gamma \notin I$.
By condition (2) of Definition \ref{M_2definition}
we know that either $I\cap M=\emptyset$ or $I\in M$.
If $I\cap M=\emptyset$ then $\sup(M\cap \gamma), \min(M\setminus \gamma)\notin I$.
Assume now that $I\in M$. Since $\gamma \notin I$ we have
that either $\gamma \leq \alpha$ or $\gamma >\beta$.
Suppose first that $\gamma \leq \alpha$.
Since $\alpha \in M$ it follows $\min(M\setminus \gamma) \leq \alpha$
and so $\min(M\setminus \gamma)\notin I$. Clearly, also $\sup(M\cap \gamma)\notin I$.
Suppose now $\gamma >\beta$. In that case, clearly, $\min(M\setminus \gamma)\notin I$.
Also, since $\beta \in M$ it follows that $\beta <\sup(M\cap \gamma)$
and so $\sup(M\cap \gamma)\notin I$.

\end{proof}

\begin{defi}\label{complete} Suppose $p\in \mathbb M_2$ and $M\in \mathcal M_p$.
We say that $p$ is $M$-{\em complete} if
\begin{enumerate}
\item $\sup (N\cap \omega_2)\in F_p$, for all $N\in \mathcal M_p$,
\item $\min (M\setminus \gamma), \sup (M\cap \gamma)\in F_p$, for all $\gamma \in F_p$.
\end{enumerate}
We say that $p$ is {\em complete} if it is $M$-complete, for all $M\in \mathcal M_p$.

\end{defi}

The following is straightforward.

\begin{fatto}\label{completion} Suppose $p \in \mathbb M_2$ and $M\in \mathcal M_p$.
Then there is an $M$-complete condition $q$ which is equivalent to $p$. We call
the least, under inclusion, such condition the $M$-{\em completion} of $p$.
\end{fatto}
\begin{proof}
First let $F^*=F_p\cup \{ \sup(N\cap \omega_2): N\in \mathcal M_p\}$.
Then let
$$
F_q=F^*\cup \{ \sup(M\cap \gamma): \gamma \in F^*\} \cup \{ \min(M\setminus \gamma): \gamma \in F^*\}.
$$
Let $A_q=A_p$ and $\mathcal M_q=\mathcal M_p$. It is straightforward
to check that $q=(F_q,A_q,\mathcal M_q)$ is a condition equivalent to $p$
and $M$-complete.
\end{proof}

\begin{remark}\label{fully-complete}
Note that in the above fact $q$ is $M$-complete for a single $M\in \mathcal M_p$.
We may not be able find $q$ which is complete, i.e. $M$-complete, for all $M\in \mathcal M_q$.
To see this, suppose there are $M,N\in \mathcal M_p$ such that
$$
\lim(M\cap N \cap \omega_2)\neq \lim (M\cap \omega_2) \cap \lim(N\cap \omega_2).
$$
Note that if $\gamma \in M\cap N$ then either $M\cap \gamma \subseteq N$ or
$N\cap \gamma \subseteq M$. Therefore, the least common limit of $M$ and $N$
which is not a limit of $M\cap N$ is above $\sup(M\cap N)$.
If $q$ is an extension of $p$ which is complete then $\sup (M\cap N)\in F_q$,
because $M\cap N \in \mathcal M_q$. Now, $\sup(M\cap N)\notin M\cap N$.
Let us assume, for concreteness, that $\sup (M\cap N)\notin M$.
We can define inductively a strictly increasing sequence $(\gamma_n)_n$
by setting $\gamma _0=\sup(M\cap N)$ and
$$
\gamma_{n+1} = \begin{cases} \min(M\setminus \gamma_n) &\mbox{if } n \mbox{ is even } \\
\min(N\setminus \gamma_n) & \mbox{if } n \mbox{ is odd. } \end{cases}
$$
Since, $q$ was assumed to be both $M$-complete and $N$-complete
we would have that $\gamma_n\in F_q$, for all $n$. This
means that $F_q$ would have to be infinite, which is a contradiction.
We do not know if such a pair of models can exist in a condition
in $\mathbb M$. Nevertheless, we will later present a variation of
$\mathbb M_2$ in which this situation does not occur and in which
the set of fully complete conditions is dense.
\end{remark}


We now come back to Lemma \ref{countable} and observe that it is valid also for $\mathbb{M}_2$.

\begin{lemma}\label{p^M}
Let $M \in \mathcal E^2$  and let $p \in  \mathbb{M}_2\cap M$.
Then  there is a new condition, which we will call $p^M$,
that is the smallest element of $\mathbb{M}_2$ extending $p$
such that $M\in \mathcal M_{p^M}$.
\end{lemma}
\begin{proof} If $M\in \mathcal E_1^2$ then simply let
$p^M=(F_p,A_p,\mathcal M_p \cup \{ M\})$.
If $M\in \mathcal E_0^2$, then, as in Lemma \ref{countable},
we let $\mathcal M_{p^M}$ be the closure of $\mathcal M_p\cup \{ M\}$ under
intersection. We also let $F_{p^M}=F_p$ and $A_{p^M}=A_p$.
We need to check that conditions (1) and (2) of Definition \ref{M_2definition}
are satisfied for $p^M$, but this is straightforward.
\end{proof}

Our next goal is to show that $\mathbb M_2$ is $\mathcal E^2$-strongly proper.
We first establish the following.

\begin{lemma}\label{generic}
Suppose $p \in \mathbb{M}_2$ and $M\in \mathcal M_p$.
Then $p$ is $(M,\mathbb M_2)$-strongly generic.
\end{lemma}

\begin{proof} We need to define, for each $r\leq p$ a
restriction $r|M \in M$ such that for every
$q \in M$ if $q\leq r|M$ then $q$ and $r$ are compatible.
So, suppose $r\leq p$. By replacing $r$ with its $M$-completion
we may assume that $r$ is $M$-complete.
We define
$$
r|M = (F_r\cap M, A_r \cap M ,\mathcal M_r \cap M).
$$
By Facts \ref{ctble-restriction} or \ref{unctble-restriction}
according to whether $M$ is countable or not we have that
$\mathcal M_r \cap M \in \mathbb M$ and therefore $r|M \in \mathbb M_2 \cap M$.
We need to show that for every $q\in M$ if $q\leq r|M$
then $q$ and $r$ are compatible.

If $M\in \mathcal E_1^2$ we already know that $\mathcal M_s= \mathcal M_q\cup \mathcal M_r$
is an $\in$-chain closed under intersection. Let $F_s=F_q\cup F_{r}$
and $A_s=A_q\cup A_r$. Finally, let $s=(F_s,A_s,\mathcal M_s)$.
It is straightforward to check that $s$ is a condition and $s\leq r,q$.

We now concentrate on the case $M\in \mathcal E_0^2$.
We define a condition $s$ as follows. We let $F_s=F_q\cup F_r$,
$A_s=A_q\cup A_r$ and
$$
\mathcal M_s = \mathcal M_q \cup \mathcal M_r \cup \{ Q\cap R : Q\in \mathcal M_q, R\in \mathcal M_r \}.
$$
We need to check that $s\in \mathbb M_2$.
By Lemma \ref{prop} we know
that $\mathcal M_s$ is an $\in$-chain closed under intersection.
Therefore we only need to check that (1) and (2) of Definition \ref{M_2definition}
are satisfied for $s$. First we check (1).

\begin{cla}\label{compatible}
$F_s \cap \bigcup A_s=\emptyset$.
\end{cla}

\begin{proof} It suffices to check that
$F_q \cap \bigcup A_r = \emptyset$ and $F_{r} \cap \bigcup A_q=\emptyset$.
Suppose first $\gamma \in F_q$ and $I \in A_r$. Since $M\in \mathcal M_r$
we have, by (2) of Definition \ref{M_2definition}, that either $I\cap M=\emptyset$ or
$I \in M$. If $I\cap M=\emptyset$ then, since $\gamma \in M$, we have
that $\gamma \notin I$. If $I\in M$ then $I\in A_r\cap M$ and,
since $q \leq r| M$, it follows that $I\in A_q$. Now, $q$ is
a condition, so $\gamma \notin I$.

Suppose now $\gamma \in F_{r}$ and $I\in A_q$. If $\gamma \in F_{r}\cap M$
then $\gamma \in F_q$. Therefore $\gamma \notin I$.
Suppose now $\gamma \in F\setminus M$. Since $r$ is
$M$-complete $\gamma^*=\min(M\setminus \gamma) \in F_r$.
Then $\gamma^* \in F_r \cap M$ and so $\gamma ^*\in F_q$.
Now, $I\in M$ and so if $\gamma \in I$ then $\gamma ^* \in I$,
which would be a contradiction. Therefore $\gamma \notin I$.

\end{proof}

We now turn to condition (2) of Definition \ref{M_2definition}.

\begin{cla}\label{Q-I}
If $Q\in \mathcal M_q$ and $I\in A_r$ then either $I\in Q$ or
$I\cap Q=\emptyset$.
\end{cla}
\begin{proof}  Since $M\in \mathcal M_r$ we have that either $I\in M$ or
$I\cap M =\emptyset$. If $I \in M$ then $I\in A_r\cap M$
and so $I\in A_q$. Since $q$ is a condition we have
that either $I\in Q$ or $I\cap Q=\emptyset$.
So, suppose $I\cap M=\emptyset$. If $Q \in \mathcal E_0^2$ then
$Q\subseteq M$ and so $Q\cap I=\emptyset$, as well.
If $Q\in \mathcal E_1^2$ then $Q\cap \omega_2$
is an initial segment of $\omega_2$, say $\gamma$.
Now, if $I\cap Q\neq \emptyset$ and $I\notin Q$ we would
have that $\gamma \in I$. Since $\gamma \in M$
this contradicts the fact that $I\cap M=\emptyset$.
\end{proof}

\begin{cla}\label{R-I}
If $R\in \mathcal M_r$ and $I\in A_q$ then either
$I\in R$ or $I\cap R=\emptyset$.
\end{cla}

\begin{proof}

Assume first that $R \in \mathcal E_1^2$. Then $R\cap \omega_2$
is an initial segment of $\omega_2$, say $\gamma$.
If $I \cap R \neq \emptyset$ and $I\notin R$ then
$\gamma \in I$. Now, since $r$ is $M$-complete we have
that $\gamma \in F_r$. If $\gamma \in M$ then $\gamma \in F_q$
and this would contradict the fact that $q$ is a condition.
If $\gamma \notin M$ let $\gamma^*=\min(M\setminus \gamma)$.
Then, again by $M$-completeness of $r$,  we have that
$\gamma^*\in F_r$. However, $\gamma^*\in M$ and therefore
 $\gamma^*\in F_q$. Since $I\in A_q$ and $q\in M$
we have that $I\in M$. If $\gamma \in I$
we would also have that $\gamma^* \in I$, which
contradicts the fact that $q$ is a condition.

We now consider the case $R\in \mathcal E_0^2$. We will
show by $\in ^*$-induction on the chain $\mathcal M_r$
that either $I\cap R=\emptyset$ or $I\in R$.
If $R \in M$ then $R\in \mathcal M_q$ so this is clear.
If $R\notin M$ then $R$ either belongs to
$[M,H(\aleph_2))_r$ or else belongs to
$[N\cap M,N)_r$, for some uncountable $N\in \mathcal M_r \cap M$.

Suppose $R\in [N\cap M,N)_r$, for some $N \in \pi_1(\mathcal M_r \cap M)$.
Since $I\in A_q$ and $N\in \mathcal M_q$ we have
that $I\in N$ or $I\cap N =\emptyset$.
On the other hand, $R\subseteq N$ so if
$I\cap N=\emptyset$ then also $I\cap R =\emptyset$.
If $I\in N$ then, since $q\in M$ and $I\in A_q$,
we have that $I \in M$ and so $I\in N\cap M$.
If there are no uncountable models in
the interval $[N\cap M,R)_r$ then $N\cap M\subseteq R$
and so $I\in R$.
If there is an uncountable model in this interval
let $S$ be the largest such model.
Now, $N\cap M\subseteq S$ and so $I\in S$ and $I\subseteq S$.
It follows that if $I\cap R\neq \emptyset$ then
also $I\cap R \cap S \neq \emptyset$.
Let $R^*=R\cap S$. Then $R^* \in \mathcal M_r$ and
$R^*$ is below $R$ in the $\in ^*$-ordering.
By the inductive assumption we would have
that $I\in R^*$ and so $I\in R$.
The case when $R\in [M,H(\aleph_2))_r$ is treated in
the same way.
\end{proof}

Finally, suppose $Q\in \mathcal M_q$, $R\in \mathcal M_r$ and
$I\in A_q\cup A_r$. Consider the relation between the model
$Q\cap R$ and $I$. If $I$ belongs to both $Q$ and $R$
then it belongs to $Q\cap R$. If $I$ is disjoint from
$Q$ or $R$ it is also disjoint from $Q\cap R$.
This completes the proof that $s$ is a condition.
Since $s\leq q,r$ it follows that $q$ and $r$ are compatible.

\end{proof}

Now, by Lemmas \ref{p^M} and \ref{generic} we have
the following.

\begin{teo}\label{M_2-proper}
The forcing $\mathbb{M}_2$ is $\mathcal E^2$-strongly proper.
Hence it is proper and preserves $\omega_2$.
\qed
\end{teo}

Suppose now $G$ is $V$-generic filter for the forcing notion $\mathbb M_2$.
We can define
$$
C_G=\bigcup \{ F_p : p \in G\} \mbox{ and } U_G=\bigcup \bigcup \{ A_p: p\in G\}.
$$
Then $C_G\cap U_G=\emptyset$. Moreover, by genericity,
$C_G\cup U_G =\omega_2$. Since $U_G$ is a union of open intervals it is
open in the order topology. Therefore, $C_G$ closed and, again by genericity,
it is unbounded in $\omega_2$.
Unfortunately, we cannot say much about the generic club $C_G$.
For reasons explained in Remark \ref{fully-complete},
we cannot even say that it does not contain infinite
subsets which are in the ground model.
In order to circumvent this problem, we now define a variation
of the forcing notion $\mathbb M_2$.
We start by some definitions.

\begin{defi}\label{lim-compatible}
Suppose $M,N\in \mathcal E^2$. We say that $M$ and $N$ are
$\lim$-compatible if
$$
\lim (M\cap N\cap \omega_2)=\lim(M\cap \omega_2)\cap \lim(N\cap \omega_2).
$$
\end{defi}

\begin{remark}\label{comp-unctble} Clearly, this conditions
is non trivial only if both $M$ and $N$ are countable.
We will abuse notation and write $\lim(M)$ for $\lim (M\cap \omega_2)$.
\end{remark}

We now define a version of the forcing notion $\mathbb M$.

\begin{defi}\label{M^*} Let $\mathbb M^*$ be the suborder of $\mathbb M$
consisting of conditions $p=\mathcal M_p$ such that any two models
in $\mathcal M_p$ are $\lim$-compatible.
\end{defi}

We have the following version of Lemma \ref{countable}.

\begin{lemma}\label{p^M2}
Let $M \in \mathcal E^2$  and let $p \in  \mathbb M^*\cap M$.
Then  there is a new condition, which we will call $p^M$,
that is the smallest element of $\mathbb{M}^*$ extending $p$
such that $M\in \mathcal M_{p^M}$.
\end{lemma}
\begin{proof} If $M\in \mathcal E_1^2$ then simply let
$p^M= M_p \cup \{ M\}$.
If $M\in \mathcal E_0^2$, then we let $\mathcal M_{p^M}$
be the closure of $\mathcal M_p\cup \{ M\}$ under
intersection. Then, thanks to Lemma \ref{countable}, we
just need to check that the models in $\M_{p^M}$ are
$\lim$-compatible. Suppose $P\in \pi_0(p)$. Then $P\in M$
and hence $P\subseteq M$. Therefore, $P$ and $M$ are lim-compatible.
Suppose now $P\in \pi_1(p)$. Then $P\cap \omega_2$ is an
initial segment of $\omega_2$, say $\gamma$.
Therefore
$$
\lim(M\cap P)=\lim(M\cap \gamma)=\lim(M)\cap (\gamma+1)=\lim(M)\cap \lim(P),
$$
and so $P$ and $M$ are lim-compatible. We also need to check
that, for any $P,Q\in \mathcal M_p$, the models $P\cap M$ and $Q\cap M$,
as well as $P\cap M$ and $Q$ are lim-compatible, but this is
straightforward.
\end{proof}

We now have a version of Lemma \ref{amalgam}.
\begin{lemma}\label{amalgam2}
Suppose $r\in \mathbb M^*$ and $M\in \mathcal M_r$. Let $q \in \mathbb M^*\cap M$
be such that $q\leq r\cap M$. Then $q$ and $r$ are compatible in $\mathbb M^*$.
\end{lemma}

\begin{proof}
If $M$ is uncountable then one can easily check that
$\mathcal M_s=\mathcal M_q\cup \mathcal M_r$ is $\in$-chain
closed under intersection and that any two models in $\mathcal M_s$
are lim-compatible.

Suppose now $M$ is countable and let
$$
\mathcal M_s = \mathcal M_q \cup \mathcal M_r \cup \{ Q\cap R : Q\in \mathcal M_q, R\in \mathcal M_r \}.
$$
Thanks to Lemma \ref{amalgam} we know that $\mathcal M_s$ is an $\in$-chain
closed under intersection. It remains to check that any two models in
$\mathcal M_s$ are lim-compatible.

\begin{cla}\label{comp2}
If $Q \in \pi_0(\M_q)$ and $R \in \pi_0(\M_r)$, then $Q$ and $R$ are lim-compatible.

\end{cla}
\begin{proof}
We show this by $\in^*$-induction on $R$.
Since $Q \in \M_q$ then  $Q\in M$ and, since $Q$ is countable, we have that $\lim(Q) \subseteq M$.
Moreover, since $R$ and $M$ are both in $\M_r$, we have that
$\lim(R \cap M) = \lim(R) \cap \lim(M)$, and so
$$
\lim(Q) \cap \lim(R) = \lim(Q) \cap \lim(R) \cap \lim(M) = \lim(Q) \cap \lim(R \cap M).
$$
Hence, without loss of generality we can assume $R$ to be $\in^*$-below $M$.
If $R\in M$ then $R\in \mathcal M_q$ and so $Q$ and $R$ are lim-compatible.
Assume now, $R\notin M$. Then by Fact \ref{ctble-restriction}
there is $N\in \pi_1(\mathcal M_r\cap M)$ such that $R\in [N\cap M,N)_r$.
We may also assume $Q$ is $\in^*$-below $N$, otherwise we could
replace $Q$ by $Q\cap N$.  Hence $Q \subseteq N \cap M$.
If there are no uncountable model in the interval $[N\cap M, R)_r$,
then $N\cap M\subseteq R$ and since $Q\in N\cap M$ we have $Q\in R$.
Therefore, $Q$ and $R$ are lim-compatible.
Otherwise, let $S$ be the $\in^*$-largest uncountable model in $[N\cap M,R)_r$.
Then $Q\in S$ and $S\cap \omega_2$ is an initial segment of $\omega_2$.
Let $R^*=R\cap S$. It follows that
 $\lim (R)\cap \lim(Q)=\lim (R^*)\cap \lim(Q)$.
By the inductive assumption we have that $\lim(R^*)\cap \lim(Q)= \lim(R^*\cap Q)$
and hence $\lim (R)\cap \lim (Q)=\lim(R \cap Q)$.
\end{proof}

Now, we need to check that any two models in $\mathcal M_s$ are lim-compatible.
So, suppose $S,S^*\in \mathcal M_s$. We may assume $S$ and $S^*$ are both
countable and of the form $S=Q\cap R$, $S^*=Q^*\cap R^*$, for $Q,Q^*\in \mathcal M_q$
and $R,R^*\in \mathcal M_r$.
Then
$$
\lim( (Q \cap R) \cap (Q^* \cap R^*)) = \lim((Q \cap Q^*) \cap (R \cap R^*))
$$
and by Claim \ref{comp2}
$$
\lim((Q \cap Q^*) \cap (R \cap R^*)) = \lim(Q \cap Q^*) \cap \lim(R \cap R^*),
$$
because $Q\cap Q^*\in \mathcal M_q$ and $R\cap R^*\in \mathcal M_r$.
Moreover, we have $=\lim(Q \cap Q^*) = \lim(Q) \cap \lim(Q^*)$
and $\lim(R \cap R^*) = \lim(R) \cap \lim(R^*)$, since the elements of $\M_q$, respectively $\M_r$,
are $\lim$-compatible. Finally, again by Claim \ref{comp2}, we have
\smallskip
$$
\lim(Q) \cap \lim(R) \cap \lim(Q^*) \cap \lim(R^*) = \lim(Q \cap R) \cap \lim(Q^* \cap R^*).
$$
\end{proof}

We now define a variation of the forcing $\mathbb M_2$ which will
have some additional properties.

\begin{defi}
Let $\mathbb{M}_2^*$ be the forcing notion whose elements are triples $p = (F_p, A_p, \mathcal{M}_p)$,
where $F_p \in [\omega_2]^{<\omega}$, $A_p$ is a finite collection of intervals of the
form $(\alpha, \beta]$, for some $\alpha, \beta <\omega_2$, $\mathcal{M}_p \in \mathbb{M}^*$, and
\begin{enumerate}
\item $F_p \cap \bigcup A_p=\emptyset$,
\item if $M \in \M_p$ and $I \in A_p$, then either $I \in M$ or $I \cap M = \varnothing$,
\end{enumerate}

\noindent The order on $\mathbb{M}_2^*$ is coordinatewise reverse inclusion, i.e.
$q \leq p$ if $F_p\subseteq F_q$,
$A_p\subseteq A_q$ and $\mathcal M_p \subseteq \mathcal M_q$.
\end{defi}

\begin{remark}\label{difference}
Note that the only difference between $\mathbb M_2^*$ and $\mathbb M_2$
is that for $p$ to be in $\mathbb M^*_2$ we require that $\mathcal M_p \in \mathbb M^*$,
i.e. the models in $\mathcal M_p$ are pairwise lim-compatible.
\end{remark}

We can now use Lemmas \ref{comp2} and \ref{amalgam2} to prove
the analogs of Lemmas \ref{p^M} and \ref{generic} for
$\mathbb M^*_2$. We then obtain the following.

\begin{teo}\label{strongly-proper2}
The forcing notions $\mathbb M^*_2$ is $\mathcal E^2$-strongly proper.
Hence, it is proper and preserves $\omega_2$.
\qed
\end{teo}

Let $G^*$ is $V$-generic filter for  $\mathbb M_2^*$.
As in the case of the forcing $\mathbb M_2$, we define

$$
C_G^*=\bigcup \{ F_p : p \in G^*\} \mbox{ and } U_G^*=\bigcup \bigcup \{ A_p: p\in G^*\}.
$$

As before $C_G^*$ is forced to be a club in $\omega_2$.
Our goal now is to show that it does
not contain any infinite subset from the ground model.
For this we will need the following lemma which
explains the reason for the requirement
of lim-compatibility for models $\mathcal M_p$,
for conditions $p$ in $\mathbb M_2^*$.

\begin{lemma}\label{complete-dense}
The set of complete conditions is dense in $\mathbb M_2^*$.
\end{lemma}

\begin{proof}
Consider a condition $p\in \mathbb M_2^*$. For each $M\in \mathcal M_p$
we consider functions $\mu_M,\sigma_M:\omega_2\rightarrow \omega_2$ defined
as follows:
$$
\mu_M(\alpha)=\min (M\setminus \alpha) \mbox { and } \sigma_M(\alpha) =\sup (M\cap \alpha).
$$

\noindent  To obtain a complete condition extending $p$ we first define:
$$
F^*_p=F_p \cup \{ \sup (M\cap \omega_2): M\in \mathcal M_p\}.
$$
We then let $\bar{F}_p$ be the closure of $F^*_p$ under the functions
$\mu_M$ and $\sigma_M$, for $M\in \mathcal M_p$.
Then $q=(\bar{F}_p,A_p,\mathcal M_p)$ will be the required  complete condition
extending $p$. The main point is to show the following.

\begin{cla}\label{complete-condition}
$\bar{F}_p$ is finite.
\end{cla}

\begin{proof}
Let $L=\bigcup \{ \lim(M): M\in \mathcal M_p\}$.
For each  $\gamma \in L$ let
$$
Y(p,\gamma) = \{ M\in \mathcal M_p : \gamma \in \lim (M)\}.
$$
and let $M(p,\gamma)=\bigcap Y(p,\gamma)$. Then, since $\mathcal M_p$ is closed
under intersection $M(p,\gamma)\in \mathcal M_p$. Since the models
in $\mathcal M_p$ are lim-compatible it follows that $\gamma \in \lim(M(p,\gamma))$.
Thus, $M(p,\gamma)$ is the least (under inclusion) model in $\mathcal M_p$
which has $\gamma$ as its limit point.
For each $\gamma \in L$ pick an ordinal $f(\gamma) \in M(p,\gamma)\cap \gamma$
above $\sup (F^*_p\cap \gamma)$ and $\sup (M\cap \gamma)$,
for all $M\in \mathcal M_p\setminus Y(p,\gamma)$.
For a limit $\gamma \in \omega_2\setminus L$ let
$$
f(\gamma)=\sup \{ \sup (M\cap \gamma) : M\in \mathcal M_p\}.
$$
Notice now that for any limit $\gamma$ and any $M\in \mathcal M_p$,
if $\xi \notin (f(\gamma),\gamma)$ then $\mu_M(\xi),\sigma_M(\xi)\notin (f(\gamma),\gamma)$.
Since $\bar{F}_p$ is the closure of $F^*_p$ under the functions
$\mu_M$ and $\sigma_M$, for $M \in \mathcal M_p$,  and
$F^*_p \cap (f(\gamma), \gamma)=\emptyset$, for all limit $\gamma$,
it follows that $\bar{F}_p \cap (f(\gamma),\gamma)=\emptyset$,
for all limit $\gamma$.
This means that $\bar{F}_p$ has no limit points and therefore is finite.
\end{proof}

\end{proof}

\begin{lemma}
Let $p \in \mathbb{M}_2^*$ be a complete condition, and let $\gamma \in \omega_2\setminus F_p$.
Then there is a condition  $q \leq p$ such that $\gamma \in I$, for some $I\in A_q$.
\end{lemma}
\begin{proof}
Without loss of generality we can assume that there is
an $M \in \M_p$ such that $\sup(M\cap \omega_2) > \gamma$,
otherwise we could let
$$
q = (F_p, A_p \cup \{(\eta, \gamma]\}, \M_p),
$$
\noindent for some $\eta < \gamma$ sufficiently large so
that $(\eta,\gamma]$ does not intersect any model in $\M_p$.

Now, since $\sup(M \cap \omega_2) \in F_p$, for every $M \in \M_p$,
the set $F_p \setminus \gamma$ is nonempty.
Let $\tau$ be $\min(F_p \setminus \gamma)$. Notice that for every
model $M \in \M_p$ either $\sup(M \cap \tau) < \gamma$, or
$\tau \in \lim(M)$, because
$$
\gamma < \sup(M \cap \tau) < \tau,
$$
would contradict the minimality of $\tau$. Moreover, if
$\sup(M \cap \tau) = \gamma$, then $\gamma$ would be in $F_p$,
contrary to the hypothesis of the lemma.

Let
$$
Y = \{M \in \M_p : \tau \in \lim(M)\}.
$$
Without loss of generality we can assume $Y \neq \emptyset$,
because otherwise we can let
$$
q=(F_p,A_p\cup \{ (\eta,\gamma]\},\mathcal M_p)
$$
for some $\eta$ sufficiently large so that
$(\eta,\gamma]$ avoids $\sup (M \cap \tau)$,
for every $M \in \M_p$.
Let $M_0=\bigcap Y$. Since $\mathcal M_p$ is closed under
intersection $M_0\in \mathcal M_p$. Moreover, since
any two models in $\mathcal M_p$ are lim-compatible
we have that $\tau \in \lim M_0$. Thus, $M_0$ is
itself in $Y$ and is contained in any member of $Y$.
Therefore,  if an interval $I$
belongs to $M_0$, then it belongs to every model in $Y$.
Let $\eta = \min(M_0 \setminus \gamma)$. Since $\tau \in \lim(M_0)$
we have $\gamma \leq \eta <\tau$. Since $\tau$
is the least element of $F_p$ above $\gamma$ it follows
that $\eta \notin F_p$.

\begin{cla}\label{eta}
$\sup (M_0 \cap \gamma) > \sup (F_p \cap \gamma)$.
\end{cla}
\begin{proof}

Suppose $\xi$ is an element of  $F_p \cap \gamma$.
Since $p$ is $M_0$-complete, we also have
$\min(M_0\setminus \xi) \in F_p$. Notice that
$\min(M_0\setminus \xi) \neq \eta$, since
$\eta \notin F_p$. Then
$$
\xi \leq \min(M_0 \setminus \xi) < \gamma,
$$
and so $\sup(M_0 \cap \gamma) > \xi$.
\end{proof}
\smallskip

Consider now some $M \in \M_p \setminus Y$.
Then $\tau \notin \lim(M)$ and, since $p$
is $M$-complete, we have that $\sup (M\cap \tau)\in F_p$.
Since $\tau$ is the least element of $F_p$ above
$\gamma$ it follows that $\sup (M\cap \tau)\in F_p\cap \gamma$.
Now, pick an element $\eta' \in M_0$ above $\sup (F_p\cap \gamma)$
and let $I=(\eta',\eta]$. It follows that $I\in M$, for all $M\in Y$
and $I \cap M=\emptyset$, for all $M\in \mathcal M_p\setminus Y$.
Therefore,
$$
q = (F_p, A_p \cup \{I\}, \M_p)
$$
is a condition stronger than $p$ and
$\gamma \in I$. Thus, $q$ is as required.
\end{proof}

\begin{coro}
If $G^*$ is a $V$-generic filter over $\mathbb{M}^*_2$,
then the generic club $C_G^*$ does not contain
any infinite subset which is in $V$.
\qed
\end{coro}

\section{Strong chains of uncountable functions}

We now consider the partial order $(\omega_1^{\omega_1},<_{\rm fin})$
of all functions from $\omega_1$ to $\omega_1$ ordered by
$f<_{\rm fin}g$ iff $\{ \xi : f(\xi)\geq g(\xi)\}$ is finite.
In \cite{Koszmider} Koszmider constructed a forcing notion
which preserves cardinals and adds an $\omega_2$ chain
in $(\omega_1^{\omega_1},<_{\rm fin})$. The construction
uses an $(\omega_1,1)$-morass which is a stationary coding
set and is quite involved. In this section we present a streamlined version
of this forcing which uses generalizes side conditions and
is based on the presentation of Mitchell \cite{Mitchell-Koszmider}.
Before that we show that Chang's conjecture implies that there
is no such chain. The argument is inspired by a proof of Shelah from \cite{Shelah-flat}.
A similar argument appears in \cite{Koszmider-sets}.

\begin{prop}\label{Chang_conj} Assume Chang's conjecture. Then there is no
chain in $(\omega_1^{\omega_1},<_{\rm fin})$ of length $\omega_2$.
\end{prop}

\begin{proof} Assume towards contradiction that
Chang's conjecture holds and $\{ f_\alpha: \alpha <\omega_2\}$
is a chain in $(\omega_1^{\omega_1},<_{\rm fin})$.
Given a function $g:I\rightarrow \omega_1$ and $\eta <\omega_1$ we
let $\min(g,\eta)$ be the function defined by:
$$
\min (g,\eta)(\zeta)=\min (g(\zeta),\eta).
$$
For each $\alpha <\omega_2$, and $\xi,\eta <\omega_1$ we define
a function $f_{\alpha}^{\xi,\eta}$ by:
$$
f_{\alpha}^{\xi,\eta}= \min (f_{\alpha}\restriction [\xi,\xi +\omega), \eta).
$$
Given  $\xi, \eta<\omega_1$, the sequence
$\{ f_{\alpha}^{\xi,\eta}: \alpha <\omega_2\}$ is $\leq_{\rm fin}$-increasing.
We define a club $C^{\xi,\eta}\subseteq \omega_2$ as follows.

\medskip

\noindent {\em Case 1:} If the sequence $\{ f_{\alpha}^{\xi,\eta}: \alpha <\omega_2\}$
eventually stabilizes under $=_{\rm fin}$ we let $C^{\xi,\eta}=\omega_2 \setminus \mu$,
where $\mu$ is least such that $f_{\nu}^{\xi,\eta}=_{\rm fin} f_{\mu}^{\xi,\eta}$,
for all $\nu \geq \mu$.

\medskip

\noindent {\em Case 2:} If the sequence $\{ f_{\alpha}^{\xi,\eta}: \alpha <\omega_2\}$
does not stabilize we let $C^{\xi,\eta}$ be a club in $\omega_2$ such that
$f_{\alpha}^{\xi,\eta} \lneqq_{\rm fin} f_{\beta}^{\xi,\eta}$, for all $\alpha,\beta\in C^{\xi,\eta}$
with $\alpha <\beta$. This means that for every such $\alpha$ and $\beta$
the set
$$
\{ n: f_{\alpha}(\xi +n) <f_{\beta}(\xi +n)\leq \eta \}
$$
is infinite.

\medskip
Let $C=\bigcap \{ C^{\xi,\eta}: \xi,\eta <\omega_1\}$.
Then $C$ is a club in $\omega_2$.
We define a coloring $c:[C]^2\rightarrow \omega_1$ by
$$
c \{ \alpha,\beta\}_<= \max \{ \xi: f_{\alpha}(\xi)\geq f_{\beta}(\xi)\}.
$$
By Chang's conjecture we can find an increasing $\omega_1$ sequence $S=\{  \alpha _{\rho}: \rho <\omega_1\}$
of elements of $C$ such that $c[[S]^2]$ is bounded in $\omega_1$.
Let $\xi =\sup (c[[S]^2]) +1$. Therefore for every $\rho <\tau <\omega_1$ we have
$$
f_{\alpha_\rho}\restriction [\xi,\omega_1) < f_{\alpha_\tau}\restriction [\xi,\omega_1).
$$
Now, let $\eta = \sup ({\rm ran}(f_{\alpha_1} \! \! \restriction \! [\xi,\xi +\omega)))$.
It follows that for every $n$:
$$
f_{\alpha_0}(\xi +n)< f_{\alpha_1}(\xi +n)\leq \eta.
$$
Since $\alpha_0,\alpha_1\in C^{\xi,\eta}$ it follows that $C^{\xi,\eta}$ was defined
using Case 2. Therefore the sequence $\{ f_{\alpha_\rho}^{\xi,\eta}: \rho <\omega_1\}$
is $\leq$-increasing and $f_{\alpha_\rho}\neq_{\rm fin}f_{\alpha_\tau}$, for all
$\rho <\tau$. For each $\rho <\omega_1$ let $n_\rho$ be the least such
that $f_{\alpha_\rho}^{\xi,\eta}(\xi +n_\rho) <f_{\alpha_{\rho +1}}^{\xi,\eta}(\xi +n_\rho)$.
Then there is a integer $n$ such that $X=\{ \rho <\omega_1 : n_\rho =n\}$ is uncountable.
It follows that the sequence $\{ f_{\alpha_{\rho}}(\xi +n): \rho \in X\}$
is strictly increasing. On the other hand it is included in $\eta$ which is countable,
a contradiction.

\end{proof}

Therefore, in order to add a strong $\omega_2$-chain in $(\omega_1^{\omega_1},<_{\rm fin})$
we need to assume that Chang's conjecture does not hold. In fact, we will assume
that there is an increasing function $g:\omega_1\rightarrow \omega_1$ such that

\begin{enumerate}
\item $g(\xi)$ is indecomposable, for all $\xi <\omega_1$,
\item ${\rm o.t.}(M\cap \omega_2) < g(\delta_M)$, for all $M\in \mathcal E^2_0$.
\end{enumerate}

It is easy to add such a function by a preliminary forcing. For instance, we can
add by countable conditions an increasing function $g$ which dominates
all the canonical functions $c_\alpha$, for $\alpha <\omega_2$, and such that
$g(\xi)$ is indecomposable, for all $\xi$.
Moreover, we may assume that $g$ is definable in the structure
$(H(\aleph_2),\in,\unlhd)$ and so it belongs to $M$, for all $M\in \mathcal E^2$.

Our plan is to add an $\omega_2$-chain $\{ f_\alpha: \alpha <\omega_2\}$
in $(\omega_1^{\omega_1},<_{\rm fin})$ below this function $g$. We can view this chain
as a single function $f:\omega_2\times \omega_1 \rightarrow \omega_1$.
We want to use conditions of the form $p=(f_p,\mathcal M_p)$,
where $f_p:A_p\times F_p\rightarrow \omega_1$ for some finite
$A_p\subseteq \omega_2$ and $F_p\subseteq \omega_1$, and $\mathcal M_p\in \mathbb M$
is a side condition. Suppose $\alpha,\beta \in A_p$ with $\alpha <\beta$,
and $M\in \mathcal \pi_0(\mathcal M_p)$. Then $M$ should localize the disagreement of
$f_\alpha$ and $f_\beta$, i.e. $p$ should force that
the finite set $\{ \xi : f_{\alpha}(\xi)\geq f_{\beta}(\xi)\}$
is contained in $M$. This means that if $\xi \in \omega_1\setminus M$ then
$p$ makes the commitment that $f_{\alpha}(\xi) <f_{\beta}(\xi)$.
 Moreover, for every $\eta \in (\alpha,\beta)\cap M$ we should have
that $f_{\alpha}(\xi)<f_{\eta}(\xi)<f_{\beta}(\xi)$.
Therefore, $p$ imposes that $f_{\beta}(\xi)\geq f_{\alpha}(\xi) + {\rm o.t.}([\alpha,\beta)\cap M)$.
This motivates the definition of the distance function below.
Before defining the distance function we need to prove some
general properties of side conditions. For a set of ordinals $X$
we let $\overline{X}$ denote the closure of $X$ in the order topology.

\begin{fatto}\label{comp1}
Suppose $P,Q\in \mathcal E^2_0$ and $\delta_P\leq \delta_Q$.
\begin{enumerate}
\item If $\gamma \in P\cap Q\cap \omega_2$
then $P\cap \gamma \subseteq Q\cap \gamma$.
\item If $P$ and $Q$ are $\lim$-compatible and
$\gamma \in \overline{P\cap \omega_2}\cap \overline{Q\cap \omega_2}$
then  $P\cap \gamma \subseteq Q\cap \gamma$.
\end{enumerate}
\end{fatto}

\begin{proof}
\noindent (1) For each $\alpha <\omega_2$ let $e_\alpha$ be the $\unlhd$-least
injection from $\alpha$ to $\omega_1$. Then $P\cap \gamma =e_\gamma^{-1}[\delta_P]$
and $Q\cap \gamma=e_\gamma^{-1}[\delta_Q]$. Since $\delta_P\leq \delta_Q$ we
have that $P\cap \gamma \subseteq Q\cap \gamma$.

\noindent (2) If $\gamma \in P\cap Q$ this is (1). Suppose
$\gamma$ is a limit point of either $P$ or $Q$ then it is also the limit point of the other.
Since $P$ and $Q$ are $\lim$-compatible we have that $\gamma \in \lim(P\cap Q)$.
Then $P\cap \gamma =\bigcup \{ e_\alpha^{-1}[\delta_P]: \alpha \in P\cap Q\}$
and $Q\cap \gamma =\bigcup \{ e_\alpha^{-1}[\delta_Q]: \alpha \in P\cap Q\}$
Since $\delta_P\leq \delta_Q$ we conclude that $P\cap \gamma \subseteq Q\cap \gamma$.
\end{proof}

\begin{fatto}\label{comp-inter}
Suppose $p\in \mathbb M$ and $P,Q\in \pi_0(\mathcal M_p)$. If $\delta_P<\delta_Q$
and $P\subseteq Q$ then $P\in Q$.
\end{fatto}
\begin{proof}
If there is no uncountable model in the interval $(P, Q)_p$, then
$P \in Q$ by transitivity. Otherwise, let
$S$ be the $\in^*$-largest uncountable model below $Q$ and we proceed by $\in^*$-induction.
First note that $S\in Q$ by transitivity and if we let $Q^*=Q\cap S$ then
$\delta_{Q^*}=\delta_Q$.
Since $P \subseteq S$, we have that $P\subseteq Q^*$ and so $Q^*$ is $\in^*$-above $P$.
By the inductive assumption we have $P\in Q^*\subseteq Q$, as desired.
\end{proof}


\begin{defi}\label{linked}
Let $p = \M_p \in \mathbb{M}^*$, $\alpha, \beta \in \omega_2$ and let
$\xi$ be a countable ordinal. Then the binary relation $L_{p,\xi}(\alpha, \beta)$ holds
if there is a $P \in \M_p$, with $\delta_P\leq \xi$, such that
$\alpha, \beta \in \overline{P\cap \omega_2}$. In this case we will say
that $\alpha$ and $\beta$ are $p,\xi$-{\em linked}.
\end{defi}

\begin{defi}\label{connected}
Let $p=\mathcal M_p\in \mathbb{M}^*$ and $\xi <\omega_1$. We let $C_{p,\xi}$ be the transitive
closure of the relation $L_{p,\xi}$. If $C_{p,\xi}(\alpha,\beta)$ holds we say
that $\alpha$ and $\beta$ are $p,\xi$-{\em connected}. If $\alpha <\beta$ and
$\alpha$ and $\beta$ are $p,\xi$-connected we write $\alpha <_{p,\xi}\beta$.
\end{defi}

From Fact~\ref{comp1}(2) we now have the following.

\begin{fatto}\label{transitive}
Suppose $p=\mathcal M_p\in \mathbb{M}^*$ and $\xi<\omega_1$.
\begin{enumerate}
\item Suppose $\alpha <\beta <\gamma$ are ordinal in $\omega_2$.  If $L_{p,\xi}(\alpha,\gamma)$ and $L_{p,\xi}(\beta,\gamma)$
hold, then so does $L_{p,\xi}(\alpha,\beta)$.
\item If $\alpha<_{p,\xi}\beta$ then there is a sequence $\alpha=\gamma_0<\gamma_1<\ldots <\gamma_{n}=\beta$
such that $L_{p,\xi}(\gamma_i,\gamma_{i+1})$ holds, for all $i<n$.
\end{enumerate}
\qed
\end{fatto}

We now present some properties of the relation $<_{p, \xi}$, in order
to define the distance function we will use in the definition of
the main forcing.

\begin{fatto}\label{tree}
Let $p=\mathcal M_p\in \mathbb{M}^*$ and $\xi <\omega_1$.
Suppose $\alpha < \beta < \gamma<\omega_2$. Then
\begin{enumerate}
\item if $\alpha <_{p,\xi} \beta$ and $\beta <_{p,\xi} \gamma$, then $\alpha <_{p,\xi}\gamma$,
\item if $\alpha <_{p,\xi}\gamma$ and $\beta <_{p,\xi} \gamma$, then $\alpha <_{p,\xi} \beta$.
\end{enumerate}
\end{fatto}
\begin{proof}
Part (1) follows directly from the definition of the relation $<_{p,\xi}$.
To prove (2) let $\alpha=\gamma_0 < \ldots < \gamma_{n}=\gamma$ witness the $p,\xi$-connection between
$\alpha$ and $\gamma$ and let $\beta=\delta_0 < \ldots < \delta_{l}=\gamma$  witness
the $p,\xi$-connection between $\beta$ and $\gamma$. We have that
$L_{p,\xi}(\gamma_i,\gamma_{i+1})$ holds, for all $i<n$, and
$L_{p,\xi}(\delta_j,\delta_{j+1})$ holds, for all $j<l$.
We prove that $\alpha$ and $\beta$
are $p,\xi$-connected by induction on $n+l$. If $n=l=1$ this is simply
Fact~\ref{transitive}(1). Let now $n,l>1$. Assume for concreteness that
$\delta_{l-1}\leq \gamma_{n-1}$. By Fact~\ref{transitive}(1)
$L_{p,\xi}(\delta_{l-1},\gamma_{n-1})$ holds; so $\alpha <_{p,\xi}\gamma_{n-1}$ and
$\beta <_{p,\xi}\gamma_{n-1}$. Now, by the inductive assumption we conclude
that $\alpha$ and $\beta$ are $p,\xi$-connected, i.e. $\alpha<_{p,\xi} \beta$.
The case $\gamma_{n-1}<\delta_{l-1}$ is treated similarly.
\end{proof}

The above lemma shows in $(1)$ that the relation $<_{p,\xi}$ is transitive and
in $(2)$ that the set $(\omega_2, <_{p, \xi})$ has a tree structure.
Since for every $M\in \mathcal E^2_0$ if
$\delta_M\leq \xi$ then ${\rm o.t.}(M\cap \omega_2) <g(\xi)$ and
$g(\xi)$ is indecomposable we conclude that the height of
$(\omega_2,<_{p,\xi})$ is at most $g(\xi)$. For every $\alpha <_{p,\xi}\beta$
we let $(\alpha,\beta)_{p,\xi}=\{ \eta: \alpha<_{p,\xi}\eta<_{p,\xi}\beta\}$.
We define similarly $[\alpha,\beta)_{p,\xi}$ and $(\alpha,\beta]_{p,\xi}$
and $[\alpha,\beta]_{p,\xi}$. If $0<_{p,\xi}\beta$, i.e. $\beta$ belongs
to some  $M\in \mathcal M_p$ with $\delta_M\leq \xi$ we write $(\beta)_{p,\xi}$ for the interval
$[0,\beta)_{p,\xi}$. Thus, $(\beta)_{p,\xi}$ is simply the set of predecessors
of $\beta$ in $<_{p,\xi}$. If $\beta$ does not belong to $\overline{M\cap \omega_2}$ for any
$M\in \mathcal M_p$ with $\delta_M\leq \xi$ we leave $(\beta)_{p,\xi}$ undefined.
 Note that when defined $(\beta)_{p,\xi}$ is a closed
subset of $\beta$ in the ordinal topology.

\begin{fatto}\label{connected}
Let $p \in \mathbb{M}^*$, $M \in \M_p$, $\xi \in M\cap \omega_1$ and
$\beta \in M \cap \omega_2$. Then $(\beta)_{p,\xi}\subseteq M$.
Moreover, if we let $p^*=p\cap M$ then $(\beta)_{p,\xi}=(\beta)_{p^*,\xi}$.
\end{fatto}
\begin{proof}
Let $\alpha <_{p,\xi}\beta$ and fix a sequence $\alpha =\gamma_0<\gamma_1<\ldots <\gamma_n=\beta$
such that $L_{p,\xi}(\gamma_i,\gamma_{i+1})$ holds, for all $i<n$.
We proceed by induction on $n$. Suppose first $n =1$ and let
$P$ witness that $\alpha$ and $\beta$ are $p,\xi$-linked.
Since $\delta_{P} <\delta_M$ we have by Fact~\ref{comp1} that $P\cap \beta \subseteq M$
and by Fact~\ref{comp-inter} that $P\cap M \in M$. Therefore
$\alpha,\beta \in \overline{P\cap M\cap \omega_2}\subseteq M$ and so $P\cap M$
witnesses that $\alpha$ and  $\beta$ are $p^*,\xi$-linked.
Consider now the case $n>1$. By the same argument as in the case $n=1$
we know that $\gamma_{n-1}$ and $\beta$ are $p^*,\xi$-linked
and then by the inductive hypothesis we conclude that $\alpha$ and
$\beta$ are $p^*,\xi$-connected.
\end{proof}

\begin{fatto}\label{factor}
Let $p \in \mathbb{M}^*$, $M \in \M_p$,
$\beta \in \omega_2 \setminus M$ and $\xi \in M \cap \omega_1$.
If $(\beta)_{p,\xi}\cap M$ is non empty then it has a largest element, say $\eta$.
Moreover, there is $Q\in \mathcal M_p\setminus M$ with $\delta_Q\leq \xi$ such
that $\eta=\sup(Q\cap M\cap \omega_2)$.
\end{fatto}
\begin{proof}
Assume $(\beta)_{p,\xi}\cap M$ is non empty and let $\eta$ be its supremum.
Note that $\eta$ is a limit ordinal. Since $(\beta)_{p,\xi}$ is a closed
subset of $\beta$ in the order topology we know that either $\eta <_{p,\xi}\beta$ or $\eta=\beta$.
By Fact~\ref{connected} $(\beta)_{p,\xi}\cap M =(\beta)_{p,\xi}\cap \eta =(\eta)_{p,\xi}$.
For every $\rho \in (\eta)_{p,\xi}$ there is some $P\in \mathcal M_p\cap M$
with $\delta_P\leq \xi$ such that $\rho \in \overline{P\cap \omega_2}$.
Since $\mathcal M_p\cap M$ is finite there is such $P$
with $\eta \in \overline{P\cap \omega_2}$. Since $P\in M$ it
follows that $\overline{P}\subseteq M$, so $\eta \in M$ and therefore
$\eta <\beta$. Finally, since $\eta$ and $\beta$ are $p,\xi$-connected, there
is a chain $\eta=\gamma_0<\gamma_1<\ldots<\gamma_n=\beta$ such that
$\gamma_i$ and $\gamma_{i+1}$ are $p,\xi$-linked, for all $i$.
Let $Q$ witness that $\eta=\gamma_0$ and $\gamma_1$ are
$p,\xi$-linked. Then $\delta_Q\leq \xi$ and $\eta =\sup(Q\cap M\cap \omega_2)$.
Since $\gamma_1\in \overline{Q\cap \omega_2}\setminus M$ it follows
that $Q\notin M$. Therefore, $Q$ is as required.
\end{proof}

We are now ready to define the distance function.

\begin{defi}
Let $p = \M_p \in \mathbb{M}^*$, $\alpha, \beta \in \omega_2$, and $\xi \in \omega_1$.
If $\alpha <_{p,\xi} \beta$ we define the $p, \xi$-distance of
$\alpha $ and $\beta$ as
$$
d_{p,\xi}(\alpha, \beta)= {\rm o.t.}([\alpha,\beta)_{p,\xi}).
$$
Otherwise we leave $d_{p,\xi}(\alpha, \beta)$ undefined.
\end{defi}

\begin{remark} Notice that for every $p$ and $\xi$ the function
$d_{p,\xi}$ is additive, i.e. if $\alpha <_{p,\xi}\beta <_{p,\xi}\gamma$ then
$$
d_{p,\xi}(\alpha,\gamma)=d_{p,\xi}(\alpha,\beta) +d_{p,\xi}(\beta,\gamma).
$$
Moreover, we have that $d_{p,\xi}(\alpha,\beta)<g(\xi)$, for every $\alpha <_{p,\xi}\beta$.
\end{remark}

We can now define the notion of forcing which adds an $\omega_2$
chain in $(\omega_1^{\omega_1},<_{\rm fin})$ below the function $g$.

\begin{defi}\label{M_3def}
Let $\mathbb{M}_3^*$ be the forcing notion whose elements are pairs $p = (f_p, \mathcal{M}_p)$, where
$f_p$ is a partial function from $\omega_2\times \omega_1$ to $\omega_1$,
$\dom (f_p)$ is of the form $A_p\times F_p$ where $0 \in A_p\in [\omega_2]^{<\omega}$,
$F_p\in [\omega_1]^{<\omega}$, $\mathcal{M}_p \in \mathbb{M}^*$, and
for every $\alpha,\beta \in A_p$ with $\alpha <\beta$, every $\xi \in F_p$ and $M\in \mathcal M_p$:
\begin{enumerate}
\item  $f_p(\alpha,\xi)<g(\xi)$,
\item if $\alpha <_{p,\xi} \beta$ then $f_p(\alpha,\xi)+ d_{p,\xi}(\alpha,\beta)\leq f_p(\beta,\xi)$,
\end{enumerate}
We let $q\leq p$ if $f_p\subseteq f_q$, $\mathcal M_p\subseteq \mathcal M_q$ and
for every $\alpha,\beta \in A_p$ and  $\xi \in F_q\setminus F_p$
if $\alpha <\beta$ then $f_q(\alpha,\xi) <f_q(\beta,\xi)$.
\end{defi}

We first show that for any $\alpha <\omega_2$ and $\xi <\omega_1$
any condition $p\in \mathbb{M}^*_3$ can be extended
to a condition $q$ such that $\alpha \in A_q$ and $\xi \in F_q$.

\begin{lemma}\label{density2}
Let $p \in \mathbb{M}_3^*$ and $\delta \in \omega_2 \setminus A_p$. Then there is
a condition $q \leq p$ such that $\delta \in A_q$.
\end{lemma}
\begin{proof}
We let $\mathcal M_q=\mathcal M_p$, $A_q=A_p\cup \{ \delta\}$ and $F_q=F_p$.
On $A_p\times F_p$ we let $f_q$ be equal to $f_p$.
We  need to define $f_q(\delta, \xi)$, for $\xi \in F_p$.
Consider one such $\xi$. If $\delta$ does not belong to $\overline{M\cap \omega_2}$,
for any $M\in \mathcal M_q$ with $\delta_M\leq \xi$, we can define $f_q(\delta,\xi)$ arbitrarily.
Otherwise, we need to ensure that if $\alpha\in A_p$ and $\alpha<_{p,\xi}\delta$
then
$$
f_{p,\xi}(\alpha,\xi) + d_{p,\xi}(\alpha,\delta)\leq f_q(\delta,\xi).
$$
Similarly, if $\beta \in A_p$ and $\delta<_{p,\xi}\beta$ we have to ensure that
$$
f_q(\delta,\xi) + d_{p,\xi}(\delta,\beta)\leq f_p(\beta,\xi).
$$
By the additivity of $d_{p,\xi}$ we know that if $\alpha <_{p,\xi}\delta <_{p,\xi}\beta$
then $d_{p,\xi}(\alpha,\beta)=d_{p,\xi}(\alpha,\delta)+d_{p,\xi}(\delta,\beta)$.
Since $p$ is a condition we know that if $\alpha, \beta \in A_p$ then
$f_p(\beta,\xi)\geq f_p(\alpha,\xi) + d_{p,\xi}(\alpha,\beta)$.
Let $\alpha^*$ be the largest element of $A_p\cap (\delta)_{p,\xi}$.
We can then simply define $f_q(\delta,\xi)$ by
$$
f_q(\delta,\xi)=f_p(\alpha^*,\xi) + d_{p,\xi}(\alpha^*,\delta).
$$
It is straightforward to check that the $q$ thus defined is a condition.
\end{proof}

\begin{lemma}\label{density3}
Let $p \in \mathbb{M}_3^*$ and $\xi \in \omega_1 \setminus F_p$. Then there is
a condition $q \leq p$ such that $\xi \in F_q$.
\end{lemma}

\begin{proof}
We let $\mathcal M_q =\mathcal M_p$, $A_q=A_p$ and $F_q=F_p\cup \{ \xi\}$.
Then we need to extend $f_p$ to $A_q\times \{ \xi\}$.
Notice that we now have the following commitments. Suppose $\alpha,\beta\in A_p$ and
$\alpha <\beta$, then we need to ensure that $f_q(\alpha,\xi) <f_q(\beta,\xi)$
in order for $q$ to be an extension of $p$. If in addition $\alpha<_{p,\xi}\beta$
then we need to ensure that
$$
f_q(\alpha,\xi)+d_{q,\xi}(\alpha,\beta)\leq  f_q(\beta,\xi)
$$
in order for $q$ to satisfy (2) of Definition~\ref{M_3def}.
We define $f_q(\beta,\xi)$ by induction on $\beta \in A_q$ as follows.
We let $f_q(0,\xi)=0$. For $\beta >0$ we let $f_q(\beta,\xi)$ be the maximum
of the following set:
$$
\{ f_q(\alpha,\xi)+1: \alpha \in (A_q \cap \beta)\setminus (\beta)_{q,\xi}\} \cup
\{ f_q(\alpha,\xi)+d_{q,\xi}(\alpha,\beta): \alpha \in A_q \cap (\beta)_{q,\xi}\}.
$$
It is easy to see that $f_q(\beta,\xi)<g(\xi)$, for all $\beta \in A_q$,
and that $q$ is a condition extending $p$.
\end{proof}

 In order to prove
strong properness of $\mathbb{M}^*_3$ we need to restrict to a relative club
subset of $\mathcal E^2$ of elementary submodels of $H(\aleph_2)$ which
are the restriction to $H(\aleph_2)$ of an elementary submodel of
$H({2^{\aleph_1}}^+)$.

\begin{defi} Let $\mathcal D^2$ be
the set of all $M \in \mathcal E^2$ such that $M=M^*\cap H(\aleph_2)$,
for some $M^*\prec H({2^{\aleph_1}}^+)$. We let $\mathcal D^2_0=\mathcal D^2\cap \mathcal E^2_0$
and $\mathcal D^2_1=\mathcal D^2\cap \mathcal E^2_1$.
\end{defi}

We split the proof that $\mathbb M^*_3$ is $\mathcal D^2$-strongly proper in two lemmas.

\begin{lemma}\label{M_3proper}
Let $p \in \mathbb{M}_3^*$ and $M \in \mathcal M_p\cap \mathcal D^2_0$. Then $p$
is an  $(M, \mathbb{M}_3^*)$-strongly generic condition.
\end{lemma}
\begin{proof}
Given $r \leq p$ we need to find a condition $r|M \in M$
such that  every $q \leq r|M$ which is in $M$ is compatible with $r$.
By Lemma~\ref{density3} we may assume that $\sup(P)\in A_r$, for every $P\in \mathcal M_r$.
The idea is to choose $r|M$ which has the same type as $r$
over some suitably chosen parameters in $M$.
Let $D=\{ \delta_P: P\in \mathcal M_r\} \cap M$.
Since $M\in \mathcal D^2_0$ there is $M^*\prec H({2^{\aleph_1}}^+)$
such that $M=M^*\cap H(\aleph_2)$. By elementary of $M^*$, we can find
in $M$ an $\in$-chain $\mathcal M_{r^*} \in \mathbb{M}^*$ extending
$\mathcal M_r\cap M$, a finite set $A_{r^*}\subseteq \omega_2$
and an order preserving bijection $\pi:A_r\rightarrow A_{r^*}$  such:
\begin{enumerate}
\item $\pi$ is the identity function on $A_r\cap M$,
\item if $\alpha,\beta \in A_r$ then, for every $\xi \in D$,
$$
d_{r^*,\xi}(\pi(\alpha),\pi(\beta))=d_{r,\xi}(\alpha,\beta).
$$
\end{enumerate}
By Lemma~\ref{density2} we can extend $f_r\restriction (A_r\cap M) \times (F_r\cap M)$
to a function $f_{r^*}:A_{r^*}\times (F_r\cap M)\rightarrow \omega_1$ such that $(f_{r^*},\mathcal M_{r^*})$
is a condition in $\mathbb M^*_3$.  Finally, we set $r|M=r^*$.

Suppose now $q \leq r|M$ and $q\in M$. We need to find a common extension $s$ of $q$ and $r$.
We  define $\M_s$ to be the closure under intersection
of $\M_r \cup \M_q$. Indeed Lemma \ref{amalgam2} shows that
$\M_s \in \mathbb{M}^*$.
We first compute the distance function
$d_{s, \xi}$ in terms of $d_{r,\xi}$ and $d_{q,\xi}$, for $\xi<\omega_1$.
First notice that the new models which are obtained by closing $\M_q \cup \M_r$
under intersection do not create new links and therefore do not influence
the computation of the distance function.

Now, consider an ordinal $\xi <\omega_1$. If $\xi \geq \delta_M$ then
all ordinals in $\overline{M\cap \omega_2}$ are pairwise $r,\xi$-linked.
The countable models of $\mathcal M_q\setminus \mathcal M_r$ are all included
in $M$ so they do not add any new $s,\xi$-links. It follows
that in this case $d_{s,\xi}=d_{r,\xi}$.
Consider now an ordinal $\xi <\delta_M$. By Fact~\ref{connected}
if $\beta \in M$ then $(\beta)_{s,\xi}=(\beta)_{q,\xi}$.
If $\beta \notin M$ then, by Fact~\ref{factor} there is a $\eta \in A_r\cap M$
such that $(\beta)_{s,\xi}\cap M=(\eta)_{q,\xi}$. Let $\xi^*=\max(D \cap (\xi+1))$.
Then, again by Fact~\ref{factor}, $\eta$ and $\beta$ are $r,\xi^*$-connected and

$$
d_{s,\xi}(\alpha,\beta)=d_{q,\xi}(\alpha,\eta)+d_{r,\xi^*}(\eta,\beta).
$$

\medskip

Let $A_s = A_q \cup A_r$ and $F_s = F_q \cup F_r$.
Our next goal is to define an extension, call it $f_s$, of $ f_q \cup f_r$ on
$A_s\times F_s$. It remains to define $f_s$ on
$$
((A_q \setminus A_r) \times (F_r\setminus F_q)) \cup ((A_r\setminus A_q) \times (F_q\setminus F_r)).
$$

\medskip

\noindent \emph{Case 1}: Consider first $\xi \in F_r\setminus F_q$ and let us
define $f_s$ on $(A_q\setminus A_r)\times \{ \xi\}$.
We already know that $d_{s,\xi}=d_{r,\xi}$, so we
need to ensure that if $\alpha,\beta \in A_s$ and $\alpha <_{s,\xi}\beta$ then
$$
f_s(\alpha,\xi)+d_{r,\xi}(\alpha,\beta)\leq f_s(\beta,\xi).
$$
Notice that all the ordinals of $A_q$ are $r,\xi$-linked as witnessed by $M$
so then we will also have that for every $\alpha,\beta\in A_q$, if $\alpha <\beta$ then
$f_s(\alpha,\xi)<f_s(\beta,\xi)$.
In order to define $f_s(\alpha,\xi)$, for $\alpha \in A_q$, let
$\alpha^*$ be the maximal element of $(\alpha)_{r,\xi}\cap A_r$ and let
$f_s(\alpha,\xi)=f_r(\alpha^*,\xi)+d_{r,\xi}(\alpha^*,\alpha)$.
It is straightforward to check that (2) of Definition~\ref{M_3def} is satisfied
in this case.

\medskip

\noindent \emph{Case 2}: Consider now some $\xi\in F_q\setminus F_r$. What we
have to arrange is that $f_s(\alpha,\xi)<f_s(\beta,\xi)$, for every
$\alpha,\beta \in A_r$ with $\alpha <\beta$. Moreover, for every
$\alpha,\beta\in A_s$ with $\alpha <_{s,\xi}\beta$ we have to arrange
that
$$
 f_s(\alpha,\xi)+d_{s,\xi}(\alpha,\beta)\leq f_s(\beta,\xi).
$$
We define $f_s$ on $(A_r\setminus A_q)\times \{ \xi\}$ by setting
$$
f_s(\beta,\xi)=f_q(\pi(\beta),\xi).
$$

First, we show that the function $\alpha \mapsto f_s(\alpha,\xi)$ is
order preserving on $A_r$. To see this observe that, since $q\leq r^*=r|M$
and $\xi\notin F_{r^*}$, the function $\alpha \mapsto f_q(\alpha,\xi)$ is
strictly order preserving on $A_{r^*}$. Moreover, $\pi$ is order preserving
and the identity on $A_r\cap M=A_{r}\cap A_q$.

Assume now $\alpha,\beta \in A_s$ and $\alpha <_{s,\xi}\beta$.
If $\alpha,\beta \in A_q$ then, since $q$ is a condition,
$f_s(\beta,\xi)\geq f_s(\alpha,\xi)+d_{q,\xi}(\alpha,\beta)$.
On the other hand, we know that $d_{s,\xi}(\alpha,\beta)=d_{q,\xi}(\alpha,\beta)$,
so we have the required inequality in this case.
By Fact~\ref{connected} the case $\alpha\in A_r\setminus A_q$ and $\beta\in A_q$ cannot happen.
Suppose $\alpha\in A_q$ and $\beta\in A_r\setminus A_q$.
Let $\xi^*=\max(D\cap (\xi +1))$. By Fact~\ref{factor} there is $\eta \in A_r\cap M$
such that
$$
d_{s,\xi}(\alpha,\beta)=d_{s,\xi}(\alpha,\eta)+d_{r,\xi^*}(\eta,\beta).
$$
By property (2) of $\pi$ we have that $d_{r^*,\xi^*}(\eta,\pi(\beta))=d_{r,\xi^*}(\eta,\beta)$.
Since $q$ extends $r^*$ it follows that $d_{q,\xi^*}(\eta,\pi(\beta))\geq d_{r^*,\xi^*}(\eta,\pi(\beta))$.
Moreover, $q$ is a condition and so:
$$
f_q(\pi(\beta),\xi) \geq f_q(\alpha,\xi) +d_{q,\xi}(\alpha,\pi(\beta))\geq  f_q(\alpha,\xi) +d_{q,\xi^*}(\alpha,\pi(\beta)).
$$
Therefore,
$$
f_q(\pi(\beta),\xi)\geq f_q(\alpha,\xi) +d_{s,\xi}(\alpha,\beta).
$$
The final case is when $\alpha,\beta \in A_r\setminus A_q$ and $\alpha <_{s,\xi}\beta$.
Note that in this case, $\alpha$ and $\beta$ are already $r,\xi$-connected, in fact,
they are $r,\xi^*$-connected, where as before $\xi^*=\max(D\cap (\xi +1))$.
By property (2) of $\pi$ we have that $\pi(\alpha)$ and $\pi(\beta)$ are $r^*,\xi^*$-connected and
$$
d_{r^*,\xi^*}(\pi(\alpha),\pi(\beta)=d_{r,\xi}(\alpha,\beta).
$$
Since $\xi^*\leq \xi$ and $q$ extends $r^*$ we have that
$$
d_{q,\xi}(\pi(\alpha),\pi(\beta)) \geq d_{r^*,\xi^*}(\pi(\alpha),\pi(\beta)).
$$
Since $q$ is a condition we have
$$
f_q(\pi(\beta),\xi)\geq f_q(\pi(\alpha),\xi) +d_{q,\xi}(\pi(\alpha),\pi(\beta)).
$$
Since $d_{s,\xi}(\pi(\alpha),\pi(\beta))=d_{q,\xi}(\pi(\alpha),\pi(\beta))$
we have $f_s(\beta,\xi)\geq f_s(\alpha,\xi) +d_{s,\xi}(\alpha,\beta)$, as required.

It follows that $s$ is a condition which extends $q$ and $r$.
This completes the proof of Lemma~\ref{M_3proper}.
\end{proof}

\begin{lemma}
Let $p \in \mathbb{M}_3^*$ and $M \in \pi_1(\mathcal M_p)$. Then $p$
is  $(M, \mathbb{M}_3^*)$-strongly generic.
\end{lemma}
\begin{proof}
Let $r \leq p$. We need to find a condition $r|M \in M$ such that
any $q \leq r|M$ in $M$ is compatible with $r$. We simply set
$$
r|M =  (f_r \up (A_r\times F_r)\cap M, \M_p \cap M).
$$
We need to show that if $q \leq r| M$ is in $M$, then there is a condition
$s \leq q, r$.
Thanks to Lemma \ref{amalgam2} we just need to define $f_s$, since we already
know that $\M_r \cup \M_q$ is an $\in$-chain and belongs to $\mathbb M^*$.
Since $\omega_1 \subseteq M$ we have that $F_r\subseteq M$ so we only
need to define an extension $f_s$ on $A_r\setminus A_q \times F_q\setminus F_r$.
We know that $M\cap \omega_2$ is an initial segment of $\omega_2$ so
all the elements of $A_r\setminus A_q=A_r\setminus M$ are above
all the ordinals of $A_q$.
Given an ordinal $\xi\in F_q\setminus F_r$ we define $f_s(\beta,\xi)$,
for $\beta \in A_r\setminus A_q$ by induction.
We set:
$$
f_s(\beta, \xi) = \max(\{f_s(\alpha, \xi)+1 : \alpha \in A_r \cap \beta\}\cup\{f_s(\alpha, \xi)+d_{s, \xi}(\alpha, 	\beta) : \alpha <_{s,\xi}\beta\}.
$$
It is easy to check that $(f_s,\mathcal M_s)$ is a condition
which extends both $q$ and $r$.
\end{proof}

\begin{coro}
The forcing $\mathbb{M}_3^*$ is $\mathcal{D}^2$-strongly
proper. Hence it preserves $\omega_1$ and $\omega_2$.
\qed
\end{coro}

We have shown that for every $\alpha <\omega_2$ and $\xi <\omega_1$ the set
$$
D_{\alpha,\xi}=\{ p\in \mathbb M ^*_3: \alpha \in A_p, \xi \in F_p\}
$$
is dense in $\mathbb M_3$. If $G$ is a $V$-generic filter in $\mathbb M^*_3$ we let
$$
f_G=\bigcup \{ f_p: p \in G\}.
$$
It follows that $f_G:\omega_2\times \omega_1 \rightarrow \omega_1$. For $\alpha <\omega_2$
we define  $f_\alpha:\omega_1\rightarrow \omega_1$ by letting
$f_\alpha(\xi)=f_G(\alpha,\xi)$, for all $\xi$. It follows that the sequence
$(f_\alpha: \alpha <\omega_2) $ is an increasing $\omega_2$-chain in
$(\omega_1^{\omega_1},<_{\rm fin})$. We have thus completed the proof of the following.

\begin{teo}There is a $\mathcal D^2$-strongly proper forcing which adds
an $\omega_2$ chain in $(\omega_1^{\omega_1},<_{\rm fin})$.
\qed
\end{teo}

\section{Thin very tall superatomic Boolean algebras}

A Boolean algebra $\mathcal B$ is called {\em superatomic }(sBa) iff every homomorphic image
of $\mathcal B$ is atomic. In particular, $\mathcal B$ is an sBa iff its Stone space
$S(\mathcal B)$ is scattered. A very useful tool for studying scattered spaces is the Cantor-
Bendixson derivative $A^{(\alpha)}$  of a set $A \subseteq S(\mathcal B)$,
defined by induction on $\alpha$ as follows.
Let $A^{(0)} = A$, $A^{(\alpha +1)}$ is the set of limit points of $A^{(\alpha)}$,
and $A^{(\lambda)}=\bigcap \{ A^{(\alpha)}: \alpha <\lambda\}$, if $\lambda$
is a limit ordinal.  Then $S(\mathcal B)$ is scattered iff for $S(\mathcal B)^{(\alpha)}=\emptyset$,
for some $\alpha$.

When this notion is transferred to the Boolean algebra $\mathcal B$, we arrive at a
sequence of ideals $I_\alpha$, which we refer to as the Cantor-Bendixson ideals, defined
by induction on $\alpha$ as follows. Let $I_0=\{ 0\}$. Given $I_\alpha$ let
$I_{\alpha +1}$ be generated by $I_\alpha$ together with all $b \in \mathcal B$
such that $b/I_\alpha$, is an atom in ${\mathcal B}/I_\alpha$.
If $\alpha$ is a limit ordinal, let $I_\alpha =\bigcup \{ I_\xi: \xi <\alpha \}$.
Then $\mathcal B$ is an sBa iff some $I_\alpha =\mathcal B$, for some $\alpha$.

The height of an sBa $\mathcal B$, ${\rm ht}(\mathcal B)$, is the least ordinal $\alpha$
such that $I_\alpha=\mathcal B$. For $\alpha <{\rm ht}{\mathcal B}$ let
${\rm wd}_{\alpha}(\mathcal B)$ be  the cardinality of the set of atoms in
${\mathcal B}/I_\alpha$.  The {\em cardinal sequence} of $\mathcal B$ is the
 sequence $({\rm wd}_\alpha(\mathcal B): \alpha <{\rm ht}(\mathcal B))$.
We say that $\mathcal  B$ is {\em $\kappa$-thin-very tall}  if ${\rm ht}(\mathcal B)=\kappa ^{++}$
 and ${\rm wd}_{\alpha}(\mathcal B)=\kappa$, for all $\alpha <\kappa^{ ++}$.
If $\kappa=\omega$ we simply say that $\mathcal B$ is {\em thin very tall}.

Baumgartner and Shelah \cite{BaumShelah} constructed a forcing notion which adds
a thin very tall sBa. This is achieved in two steps. First
they adjoin by a $\sigma$-closed $\aleph_2$-cc forcing
a function $f:[\omega_2]^2\rightarrow [\omega_2]^{\leq \omega}$ with
some special properties. Such a function is called a $\Delta$-function.
In the second step they use a $\Delta$-function to define a ccc forcing
notion which adds a thin very tall sBa. The purpose of this section is
to show how this can be achieved directly by using generalizes side conditions.
The following concept from \cite{BaumShelah} was made explicit by Bagaria in
\cite{bagaria2002lgb}.

\begin{defi}\label{theta-poset}
Given a cardinal sequence $\theta = \langle \kappa_{\alpha} : \alpha
<\lambda \rangle$, where each $\kappa_{\alpha}$ is an infinite
cardinal, we say that a structure $(T,\leq,i)$ is a $\theta$-{\em poset}
if $<$ is a partial ordering on $T$  and the following hold:
\begin{enumerate}
\item  $T = \bigcup \{ T_{\alpha}: \alpha <\lambda \}$,
where each $T_{\alpha}$ is of the form $\{\alpha\}\times
Y_{\alpha}$, and $Y_{\alpha}$ is a set of cardinality
$\kappa_{\alpha}$.
\item  If $s \in T_{\alpha}$, $t \in T_{\beta}$ and $s < t$, then $\alpha <\beta$.
\item For every $\alpha <\beta <\lambda$, if $t\in T_{\beta}$
then the set $\{ s\in T_{\alpha} : s< t\}$ is infinite.
\item  $i$ is a function from $[T]^2$ to $[T]^{<\omega}$ with
the following properties:
\begin{enumerate}
\item If $u \in i\{s, t\}$, then $u \leq s, t$
\item If $u \leq s, t$, then there exists $v \in i\{s, t\}$ such that $u \leq v$.
\end{enumerate}
\end{enumerate}
\end{defi}

We let $\Omega(\lambda)$ denote the sequence of length $\lambda$  with
all entries equal to $\omega$.
The following is implicitly due to Baumgartner (see \cite{bagaria2002lgb} for a proof).

\begin{fatto}\label{proplc}
Let $\theta= \langle \kappa_\alpha : \alpha<\lambda \rangle $ be a
sequence of cardinals. If there exists a $\theta$-poset, then there
exists an sBa whose cardinal sequence is $\theta$.
\qed
\end{fatto}

We now define a forcing notion which adds an $\Omega(\omega_2)$-poset.
If $x\in \omega_2\times \omega$ is of the form $(\alpha,n)$ then
we denote $\alpha$ by $\alpha_x$ and $n$ by $n_x$.

\begin{defi}\label{M_4def}
Let $\mathbb{M}_4$ be the forcing notion whose elements are tuples $p = (x_p,\leq_p,i_p, \mathcal{M}_p)$, where
$x_p$ is a finite subset of $\omega_2\times \omega$, $\leq_p$ is a partial
ordering on $x_p$, $i_p:[x_p]^2\rightarrow [x_p]^{<\omega}$, $\mathcal M_p \in \mathbb M$
and the following hold:
\begin{enumerate}
\item if $s,t \in x_p$ and $s <_p t$ then $\alpha_s < \alpha_t$,
\item if $u\in i_p\{ s,t\}$ then $u\leq_p s,t$,
\item for every $u\leq_p s,t$ there is $v\in i_p\{ s,t\}$ such that
$u\leq_p v$,
\item for every $s,t\in x_p$ and $M\in \mathcal M_p$ if
$s,t\in M$ then $i_p\{ s,t\}\in M$.
\end{enumerate}
We let $q\leq p$ if and only if $x_q\supseteq x_p$,
$\leq_q\restriction x_p=\leq_p$, $i_q\restriction [x_p]^2= i_p$
and $\mathcal M_p\subseteq \mathcal M_q$.
\end{defi}

We first observe that a version of Lemma \ref{countable}
holds for $\mathbb M_4$.

\begin{lemma}\label{p^M_4}
Let $M \in \mathcal E^2$  and let $p \in  \mathbb{M}_4\cap M$.
Then  there is a new condition, which we will call $p^M$,
that is the smallest element of $\mathbb{M}_4$ extending $p$
such that $M\in \mathcal M_{p^M}$.
\end{lemma}
\begin{proof} If $M\in \mathcal E_1^2$ then simply let
$p^M=(x_p,\leq_p,i_p,\mathcal M_p \cup \{ M\})$.
If $M\in \mathcal E_0^2$, then, as in Lemma \ref{countable},
we let $\mathcal M_{p^M}$ be the closure of $\mathcal M_p\cup \{ M\}$ under
intersection and let $p^M=(x_p,\leq_p,i_p,\mathcal M_{p^M})$.
We need to check that condition (4) of Definition \ref{M_4def}
is satisfied. Since $p\in M$ we have that $x_p\subseteq M$.
In the case $M\in \mathcal E_1^2$ the only new model in $\mathcal M_{p^M}$
is $M$ so condition (4) holds for $p^M$ since it holds for $p$.
In the case $M\in \mathcal E_0^2$ there are also models of the form $N\cap M$,
where $N\in \pi_1(\mathcal M_p)$. However, condition (4)
holds for both $N$ and $M$ and so it holds for their intersection.
\end{proof}

Next, we show that $\mathbb M_4$ is $\mathcal E^2$-proper. We split this in
two parts.

\begin{lemma}\label{M_4-proper-E_0}
$\mathbb M_4$ is $\mathcal E_0^2$-proper.
\end{lemma}

\begin{proof}
Let $\theta$ be a sufficiently large regular
cardinal and let $M^*$ be a countable elementary
submodel of  $H(\theta)$ containing all the relevant
objects. Then $M=M^*\cap H(\omega_2)$ belongs to $\mathcal E_0^2$.
Suppose $p\in \mathbb M_4 \cap M$.
Let $p^{M}$ be the condition defined in Lemma \ref{p^M_4},
i.e. $p^{M}=(x_p,\leq_p,i_p,\mathcal M_{p^{M}})$,
where $\mathcal M_{p^{M}}$ is the closure of
$\mathcal M_p \cup \{ M\}$ under intersection.
We show that $p^{M}$ is $(M^*,\mathbb M_4)$-generic.
Let $D \in M^*$ be a dense subset of $\mathbb M_4$ and $r\leq p^{M}$.
We need to find a condition $q\in D\cap M^*$ which is compatible
with $r$. Note that we may assume that $r\in D$.
We define a condition $r|M$ as follows.
First let $x_{r|M}=x_r\cap M$ and then let $\leq_{r|M}=\leq_r\restriction x_{r|M}$
and $i_{r|M}=i_r\restriction [x_{r|M}]^2$.
Condition (4) of Definition \ref{M_4def} guarantees that
if $s,t\in x_{r|M}$ then $i_r\{ s,t\}\subseteq M$.
Finally, let $\mathcal M_{r|M}=\mathcal M_r\cap M$.
It follows that $r|M=(x_{r|M},i_{r|M},i_{r|M},\mathcal M_{r|M})$
belongs to $\mathbb M_4\cap M$.
By elementarity of $M^*$ in $H(\theta)$ there is a condition $q\in D\cap M^*$
extending $r|M$ such that $(x_q\setminus x_{r|M})\cap N=\emptyset$,
for all $N\in \pi_0(\mathcal M_{r|M})$.

\begin{cla}\label{M_4-compatible-ctble} $q$ and $r$ are compatible.
\end{cla}
\begin{proof} We define a condition $s$ as follows.
We set $x_s=x_q\cup x_r$ and we let
$\leq_s$  be the transitive closure of $\leq_q \cup \leq_r$, i.e.
if $u\in x_q\setminus x_r$, $v\in x_r\setminus x_q$ and  $t \in x_{r|M}$ are such that
$u \leq_q t$ and $t \leq_r v$, then we let $u \leq_s v$.
Similarly, if $v\leq_r t$ and $t\leq_q u$ we let $v\leq_s u$.
We let $\mathcal M_s$ be the closure under intersection of $\mathcal M_q\cup \mathcal M_r$.
It remains to define $i_s$. For $z \in x_r$ let
$A_{z}= \{t \in x_{r|M} : t \leq_r z\}$ and for $z\in x_q$ let
$B_{z}= \{t \in x_{r|M} : t \leq_q z\}$.
We let

$$
i_s\{ u,v\}= \left\lbrace \begin{array}{ll}
                        i_q\{ u,v\} & \textnormal{if
                        $u,v\in x_q$}, \\
                        i_r\{ u,v\} & \textnormal{if
                        $u,v\in x_r$}, \\
                        \bigcup_{t \in A_v}i_q\{u,t\} \cup \bigcup_{t \in B_u}i_r\{t,v\}&
                        \textnormal{if $ u\in x_q\setminus x_r$ and $v\in x_r\setminus x_q$}. \\
 \end{array}\right.
$$

\medskip

We now need to check property (4) of Definition \ref{M_4def}, i.e. for every
$u,v\in x_s$ and $P\in \mathcal M_s$, if $u,v\in P$ then $i_s\{ u,v\} \in P$.
 First of all notice that we only need to show the above property for $P\in \M_q \cup \M_r$,
because the other models in $\M_s$ are obtained by intersection and,
if (4) holds for $u,v$ and $P$ and also for $u,v$ and $Q$, it also
holds for $u,v$ and $P\cap Q$.

\medskip

\noindent \emph{Case 1}: $u,v \in x_q$ and $P\in \mathcal M_r$.
If $P\in \mathcal M_{r|M}$ then $P\in \mathcal M_q$ and then (4)
holds since $q$ is a condition and $i_s\{ u,v\}=i_q\{u,v\}$.
Suppose now $P\in \mathcal M_r\setminus M$ and $\delta_P<\delta_M$.
Then, by Fact \ref{comp-inter}, $P\cap M\in \mathcal M_{r|M}$
and so $P\cap M\in \mathcal M_q$, therefore (4) of Definition \ref{M_4def} holds again.
Finally, if $\delta_P\geq \delta_M$ then, by Fact \ref{comp1},
$P\cap M\cap \omega_2$ is an initial segment of $M\cap \omega_2$.
We know that $i_q\{ u,v\}\in M$ and for every $w\in i_q\{ u,v\}$
$\alpha_w\leq \min (\alpha_u,\alpha_v)$. Therefore,
we have that $i_q\{u,v\}\in P\cap M$. Since $i_s\{ u,v\}=i_q\{u,v\}$,
we conclude that (4) holds in this case.

\medskip

\noindent \emph{Case 2}: $u,v \in x_r$ and $P\in \mathcal M_q$.
If $u,v \in x_{r|M}$ then $u,v\in x_q$ and
again, since $q$ is a condition and $i_s\{ u,v\}=i_q\{ u,v\}$,
we know that $i_s\{u,v\}\in P$. Suppose now that
$u$ and $v$ are not both in $M$. If $P\in \mathcal E^2_0$  then
$P\subseteq M$ and so we cannot have $u,v\in P$.
If $P\in \mathcal E^2_1$  we know that $P\cap \omega_2$ is
an initial segment of $\omega_2$. Moreover, if $w\in i_s\{u,v\}$
then $\alpha_w\leq \min(\alpha_u,\alpha_v)$ and so
if $u,v\in P$ we also have that $i_s\{u,v\}\in P$,
so (4) of Definition \ref{M_4def} holds again.

\medskip

\noindent \emph{Case 3}: $u\in x_q\setminus x_r$ and $v\in x_r\setminus x_q$.
If $P\in \mathcal E^2_1$ then $P\cap \omega_2$ is an initial segment
of $\omega_2$. Moreover, as before, we have that
$\alpha_w\leq \min(\alpha_u,\alpha_v)$, for every
$w\in i_s\{ u,v\}$. Therefore, $i_s\{ u,v\}\in P$.
Suppose now   $P\in \mathcal E^2_0$.
If $P\in \pi_0(\mathcal M_q)$ then $P\subseteq M$ so $v\notin P$.
Now assume $P\in \pi_0(\mathcal M_r)$.
If $\delta_P <\delta_M$ then by Fact \ref{comp-inter},
$P\cap M\in \mathcal M_{r|M}$. However, the condition $q$ is chosen
so that $(x_q\setminus x_{r|M})\cap N$, for all $N\in \pi_0(\mathcal M_{r|M})$,
therefore in this case $u\notin P$.
Assume now $\delta_P\geq \delta_M$. Then, by Fact \ref{comp1},
we have that $P\cap M\cap \omega_2$ is an initial segment
of $M\cap \omega_2$.
Consider first some $t\in A_v$. Then $t\in M$ and $i_q\{ u,t\}\in M$.
Moreover, $\alpha_w\leq \min(\alpha_u,\alpha_t)$, for every $w\in i_q\{ u,t\}$.
Since $P\cap M\cap \omega_2$ is an initial segment of
$M\cap \omega_2$, it follows that $\alpha_w\in P\cap M$, for
every $w\in i_q\{ u,t\}$. This implies that $i_q\{ u,t\}\in P\cap M$.
Finally, consider some $t\in B_u$. Then $t\in M$ and, since  $u\in P\cap M$,
$\alpha_t\leq \alpha_u$ and $P\cap M\cap \omega_2$ is an initial segment of $M\cap \omega_2$,
we have that $\alpha_t\in P\cap M$ and so $t\in P\cap M$.
Now, since $r$ is a condition, $P \in \mathcal M_r$ and $t,v\in P$,
we have that $i_r\{t,v\}\in P$, so (4) of Definition \ref{M_4def}
holds in this case as well.

It follows that $s$ is a condition which extends both $q$ and $r$.
This completes the proof of Claim \ref{M_4-compatible-ctble} and Lemma \ref{M_4-proper-E_0}.
\end{proof}
\end{proof}

\begin{lemma}\label{M_4-proper-E_1}
$\mathbb M_4$ is $\mathcal E_1^2$-proper.
\end{lemma}

\begin{proof}
Let $\theta$ be a sufficiently large regular
cardinal and $M^*$ an elementary submodel of  $H(\theta)$ containing
all the relevant objects such that $M=M^*\cap H(\omega_2)$
belongs to $\mathcal E_1^2$. Fix $p \in M \cap \mathbb{M}_4$.
Let $p^M$ be as in Lemma \ref{p^M_4}. We claim that $p^M$
is $(M^*,\mathbb{M}_4)$-generic. In order to verify this
consider a dense subset $D$ of $\mathbb M_4$ which belongs to
$M^*$ and a condition $r\leq p^{M}$. We need to find
a condition $q\in D\cap M^*$ which is compatible with $r$.
By extending $r$ if necessary we may assume it belongs to $D$.
Let  $r|M=(x_{r|M},i_{r|M},i_{r|M},\mathcal M_{r|M})$ be as in Lemma \ref{M_4-proper-E_0}.
By elementarity of $M^*$ in $H(\theta)$, we can find $q \leq r|M$, in $D \cap M$,
such that  $(x_q\setminus x_{r|M})\cap N=\emptyset$,
for all $N\in \mathcal M_{r|M}$, and
 if $u\in x_{r|M}$ and $v\in x_q\setminus x_{r|M}$ then $\alpha_u < \alpha_v$.

\begin{cla}\label{M_4-compatible-unctble} $q$ and $r$ are compatible.
\end{cla}

\begin{proof} We define a condition $s$ as follows.
We set $x_s=x_q\cup x_r$, $\leq_s=\leq_q\cup \leq_r$ and
$\mathcal M_s =\mathcal M_q \cup \mathcal M_r$.
Note that $\leq_s$ is a partial order and $\mathcal M_s \in \mathbb{M}$.
It remains to define $i_s$. We let

$$
i_s\{ u,v\}= \left\lbrace \begin{array}{ll}
                        i_q\{ u,v\} & \textnormal{if
                        $u,v\in x_q$}, \\
                        i_r\{ u,v\} & \textnormal{if
                        $u,v\in x_r$}, \\
                        \{ z \in x_{r|M}: z\leq_q u \mbox{ and } z \leq_r v\}&
                        \textnormal{if $ u\in x_q\setminus x_r$, $v\in x_r\setminus x_q$}. \\
 \end{array}\right.
$$
We need to check (4) of Definition \ref{M_4def}. So, suppose $u,v\in x_s$,
$P\in \mathcal M_s$ and $u,v\in P$. We need to show
that $i_s\{ u,v\} \in P$.

\medskip

\noindent \emph{Case 1}: $u,v \in x_q$ and $P\in \mathcal M_r$.
If $P\cap M\in \mathcal M_{r|M}$
this follows from the fact that $q$ is a condition and
$\mathcal M_{r|M}\subseteq \mathcal M_q$.
If $P\cap M\notin \mathcal M_{r|M}$ then
$M\subseteq P$ and, since $i_q\{u,v\}\in M$, it follows
that $i_s\{u,v\}\in P$.

\medskip

\noindent \emph{Case 2}: $u,v \in x_r$ and $P\in \mathcal M_q$.
If $u,v\in M$ then $u,v\in x_q$, so $i_s\{u,v\}\in P$ follows
from the fact that $q$ is a condition.
Assume now, for concreteness, that $v\notin M$.
If $P\in \mathcal M_q$ then $v\notin P$ and if $P\in \mathcal M_r$ then
(4) of Definition \ref{M_4def} follows from the fact that $r$ is a condition.

\medskip

\noindent \emph{Case 3}: $u\in x_q \setminus x_r$ and $v\in x_r\setminus x_q$.
If $P\in \mathcal M_q$ then $v\notin P$. If $P\in \mathcal M_r$ then
either  $P\cap M\in \mathcal M_{r|M}$ and then, by the choice of $q$, we have
that $u\notin P$. Otherwise $M\subseteq P$ and in this case
$x_{r|M}\subseteq P$. Since $i_s\{u,v\}\subseteq x_{r|M}$ it follows
that (4) of Definition \ref{M_4def} holds in this case as well.
\end{proof}
This completes the proof of  Lemma \ref{M_4-proper-E_1}.
\end{proof}

\begin{coro}
The forcing $\mathbb{M}_4$ is $\mathcal{E}^2$-proper.
Hence it preserves $\omega_1$ and $\omega_2$.
\qed
\end{coro}

It is easy to see that the set
$$
D_{\alpha, n}= \{ p\in \mathbb M_4: (\alpha,n) \in x_p \}
$$
is dense in $\mathbb{M}_4$, for every $\alpha \in \omega_2$ and $n \in \omega$.
Moreover, given $t\in \omega_2\times \omega$, $\eta <\alpha_t$ and $n <\omega$,
one verifies easily that the set
$$
E_{t,\eta,n}= \{ p: t\in x_p \mbox{ and }  |\{ i : (\eta,i)\in x_p \mbox{ and } (\eta,i)\leq_p t\}| \geq n\}
$$
is dense. Then if $G$ is $V$-generic filter on $\mathbb{M}_4$ let
$$
\leq_G = \bigcup \{ \leq_p : p \in G\} \hspace{1cm} \mbox{and} \hspace{1cm} i_G=\bigcup \{ i_p: p\in G\}.
$$
It follows that $(\omega_2\times \omega,\leq_G,i_G)$
is an $\Omega(\omega_2)$-poset in $V[G]$. We have therefore proved
the following.

\begin{teo} There is an $\mathcal E^2$-proper forcing notion
which adds an $\Omega(\omega_2)$-poset.
\qed
\end{teo}

\bibliographystyle{plain}
\bibliography{side}

\begin{thebibliography}{10}

\bibitem{bagaria2002lgb}
J.~Bagaria.
\newblock {Locally-generic Boolean algebras and cardinal sequences}.
\newblock {\em Algebra Universalis}, 47(3):283--302, 2002.

\bibitem{BaumShelah}
James~E. Baumgartner and Saharon Shelah.
\newblock Remarks on superatomic {B}oolean algebras.
\newblock {\em Ann. Pure Appl. Logic}, 33(2):109--129, 1987.

\bibitem{Fried}
Sy~David Friedman.
\newblock Forcing with side conditions.
\newblock Preprint, 2004.

\bibitem{Koszmider-sets}
Piotr Koszmider.
\newblock On the existence of strong chains in {${\mathcal P}(\omega_1)/{\rm
  Fin}$}.
\newblock {\em J. Symbolic Logic}, 63(3):1055--1062, 1998.

\bibitem{Koszmider}
Piotr Koszmider.
\newblock On strong chains of uncountable functions.
\newblock {\em Israel J. Math.}, 118:289--315, 2000.

\bibitem{Mitchell-Koszmider}
William~J. Mitchell.
\newblock Notes on a proof of koszmider.
\newblock Preprint, http://www.math.ufl.edu/~wjm/papers/koszmider.pdf, 2003.

\bibitem{Mitch}
William~J. Mitchell.
\newblock Adding closed unbounded subsets of {$\omega_2$} with finite forcing.
\newblock {\em Notre Dame J. Formal Logic}, 46(3):357--371, 2005.

\bibitem{Mitch2}
William~J. Mitchell.
\newblock {$I[\omega_2]$} can be the nonstationary ideal on {${\rm
  Cof}(\omega_1)$}.
\newblock {\em Trans. Amer. Math. Soc.}, 361(2):561--601, 2009.

\bibitem{Neeman}
Itay Neeman.
\newblock Forcing with side conditions.
\newblock slides, http://www.math.ucla.edu/~ineeman/fwsc.pdf/, 2011.

\bibitem{Shelah-flat}
Saharon Shelah.
\newblock On long increasing chains modulo flat ideals.
\newblock {\em Math. Logic Quarterly}, 56(4):397--399, 2010.

\end{thebibliography}

\end{document}